\documentclass[reqno,11pt]{article}

\usepackage{amsmath,amssymb,amsfonts}

\newcommand{\gD}{\Delta}
\newcommand{\tho}{\th_{_1}}

\newcommand{\dd}{\gD}

\newcommand{\kk}{2}

\newcommand{\sfrac}[2]{\mbox{$\frac{#1}{#2}$}}

\newcommand{\gL}{\Lambda}

\newcommand{\bb}{\Big}
\newcommand{\beq}[1]{\begin{equation}\label{#1}}
\newcommand{\mn}[0]{\medskip\noindent}
\newtheorem{thm}{Theorem}[section]
\newtheorem{cor}[thm]{Corollary}
\newtheorem{lem}[thm]{Lemma}
\numberwithin{equation}{section}
\newcommand{\gd}{\delta}
\renewcommand{\th}{\theta}
\newcommand{\case}[4]{
\left\{ \begin{array}{ll} {#1} & \mbox{#2} \\ {#3}  & \mbox{#4}
\end{array} \right.}
\newcommand{\caseth}[6]{
\left\{ \begin{array}{ll} {#1} & \mbox{#2} \\ {#3}  & \mbox{#4} \\
{#5} & \mbox{#6}  \end{array} \right.}

\newcommand{\bn}{\bigskip\noindent}
\newcommand{\sm}{\setminus}
\newcommand{\ra}{\rightarrow}
\newcommand{\pr}{\Pr}
\newcommand{\sub}{\subseteq}
\newcommand{\qed}{\hfill$\square$}
\newcommand{\gs}{\sigma}
\newcommand{\raf}[1]{(\ref{#1})}
\newcommand{\pf}{{\bf Proof. }}
\newcommand{\eps}{\varepsilon}
\newcommand{\ga}{\alpha}
\newcommand{\gb}{\beta}
\newcommand{\gl}{\lambda}

\newcommand{\bean}{\begin{eqnarray*}}
\newcommand{\eean}{\end{eqnarray*}}

\newcommand{\gt}{\theta}

\newcommand{\poi}{{\rm Poi}}
\newcommand{\eeq}{\end{equation}}
\newcommand{\old}[1]{}

\newcommand{\f}{{\cal F}}

\newcommand{\thl}{\th_{_\gl}}

\newcommand{\odvi}{\topmargin -70pt \advance \topmargin by
-\headheight \advance \topmargin by -\headsep \textheight 7.7in
\oddsidemargin -67pt \evensidemargin \oddsidemargin
\marginparwidth 0.5in \textwidth 6.5in}

\newcommand{\dvi}{}

\newcommand{\reg}{\topmargin 0pt \advance \topmargin by -\headheight \advance
\topmargin by -\headsep \textheight 9.1in \oddsidemargin -.1in
\evensidemargin \oddsidemargin \marginparwidth 0.5in \textwidth
6.5in}

\reg

\newcommand{\fnp}{F_{_{\!{PC}}}(n,p\,;2)}
\newcommand{\fnpk}{F_{_{\!{PC}}}(n,p\,;k)}

\newcommand{\npk}{F(n,p\,;k)}

\begin{document}

\vspace{6mm}
\begin{center}
{\Large \bf  Finding cores of random $2$-SAT formulae via  Poisson
cloning  } \\[10mm]
{\it Jeong Han Kim\footnote{This work was partially supported by
Yonsei University Research Funds 2006-1-0078 and 2007-1-0025, and by
the second stage of the Brain Korea 21 Project in 2007,
and by the Korea Research Foundation Grant funded by the Korean Government (MOEHRD) (KRF-2006-312-C00455).}\\
 Department of Mathematics and\\
 Random Graph Research Center\\
 Seoul, 120-749 Korea\\
 {\tt jehkim@yonsei.ac.kr}}
\end{center}
\renewcommand{\include}{\old}
\newcommand{\edc}{}
\include{dstart}


\vspace{5.5mm}

\mn{\bf Abstract.} For the random $2$-SAT formula $F(n,p)$, let
$F_C (n,p)$ be the  formula left after the pure literal algorithm
applied to $F(n,p)$ stops. Using the recently developed Poisson
cloning model together with the cut-off line algorithm (COLA), we
completely analyze the structure of $F_{C} (n,p)$. In particular,
it is shown that, for $\gl:= p(2n-1) = 1+\gs $ with $\gs\gg
n^{-1/3}$, the core of $F(n,p)$ has $\thl^2 n +O((\thl n)^{1/2})$
variables and $\thl^2 \gl n+O((\thl n))^{1/2}$ clauses, with high
probability, where $\thl$ is the larger solution  of the equation
$\th- (1-e^{-\thl \gl})=0$. We also estimate the probability of
$F(n,p)$ being satisfiable to obtain
$$
\pr[ F_2(n, \sfrac{\gl}{2n-1}) ~\mbox{is satisfiable} ] =
\caseth{1-\frac{1+o(1)}{16\gs^3 n} }{if $\gl= 1-\gs$ with $\gs\gg
n^{-1/3}$}{}{}{ e^{-\Theta(\gs^3n)}}{if $\gl=1+\gs$ with $\gs\gg
n^{-1/3}$,}
$$
where $o(1)$ goes to $0$ as $\gs$ goes to $0$. This improves the
bounds of Bollob\'as et al. \cite{BBCKW}.

\date{today}

\mn

\section{Introduction}
\label{intro}

An instance of the satisfiability problem is given by a
conjunctive normal form (CNF), that is, a conjunction of
disjunctions. Each disjunction, or clause,  is of the form
$(y_1\vee \cdots \vee y_k)$, where $y_i$'s are chosen among $2n$
literals consisting of $n$ Boolean variables, conditioned that all
$k$ literals are strictly distinct, i.e.,  no literals with the
same underlying variables appear more than once. The problem is
whether a given formula has an assignment of truth values ($0$ or
$1$) for the $n$ variables that satisfies the formula. When such
an assignment exists, the formula is called satisfiable. It is
unsatisfiable, otherwise. It is now well-known that the
satisfiability problem is NP-complete (\cite{Coo71}). Even the
$k$-satisfiability problem, in which each clause consists of
exactly $k$ literals, is known to be NP-complete for $k\geq 3$
(\cite{Coo71}). In case of  $k=2$, there is a polynomial time
algorithm \cite{Coo71} to determine wether the instance of the
$2$-satisfiability problem is satisfiable or not.

The random $k$-SAT   formula $\npk$ on $n$ variables  is  the
conjunction of clauses
 selected with probability $p$ from
 the set of $2^{k} {n \choose k} $ possible clauses,
 independent of all others.
Not surprisingly, the random $2$-SAT and the random $3$-SAT
formulae have been most extensively studied  and many research
papers regarding the random models  have been published. For
$k=2$, Chv\'{a}tal and Reed \cite{CR}, Goerdt \cite{Go} and
Fernandez de la Vega \cite{FV} independently proved that the
random $2$-SAT problem undergoes  a phase transition at $1$, that
is,
$$\lim_{n\ra \infty}
\pr[ F_2(n,\sfrac{\gl}{2n-1}) ~\mbox{is satisfiable} ] =
\case{1}{if $ \gl  <1$}{ 0}{if $\gl >1$.} $$ Though there is no
essential difference, we prefer $\gl=p(2n-1)$ to  $\gl=2pn$
because $p(2n-1)$ is the mean average degree of each literal.
Techniques used to prove the phase transition are essentially
based upon the first and the second moment methods for the number
of certain structures closely related to the satisfiability.
Bollob\'as et al. \cite{BBCKW} took much more sophisticated
approaches to determine the scaling window for the problem:
$$
\pr[ F_2(n, \sfrac{\gl}{2n-1}) ~\mbox{is satisfiable} ] =
\caseth{1-\Theta(\frac{1}{\gs^3 n}) }{if $\gl= 1-\gs$ withn
$\gs\gg n^{-1/3}$}{}{}{ e^{-\Theta(\gs^3n)}}{if $\gl=1+\gs$ with
$\gs\gg n^{-1/3}$.}
$$

Though it is believed that the random $k$-SAT problem, $k\geq 3$,
undergoes a similar phase transition, it remains as a conjecture.
Only sharp transitions are known due to a seminal result of
Friedgut \cite{FB}. The upper and lower bounds for the critical
value $\gl_3$, (assuming the conjecture is true for $k=3$) have
colorful history. In a series of papers
\cite{BFU,EF,KM,DB,Kea,Jea,Zi,Ka,DBM}, the upper bound of $\gl_3$
has been improved to  4.506. There has been considerable work
bounding $\gl_3$ from below too. The easiest but fundamental
algorithm is the pure literal algorithm (PLA).
 A literal
is {\em pure} in a formula if it belongs to at least one  clause
of the formula, while its negation is in no clause. The PLA keeps
selecting a pure literal, setting it  true, and removing clauses
containing the literal as they are already satisfied. This
procedure may (or may not) yield new pure literals. The algorithm
stops when no more pure literal is left.  We say that the PLA {\em
succeeds} if no clause remains in the formula after it stops.
Clearly, the formula is satisfiable if the PLA succeeds. The
converse is not true, for example, $(y\vee z) \wedge (\bar{y},
\bar{z})$ is satisfiable whereas no pure literal exists.

Broder, Frieze, and Upfal \cite{BFU} analyzed the PLA for the
random $3$-SAT problem to show that,   if $\gl < 1.225$ then the
PLA applied to $F(n,\frac{\gl}{n^2}\,; 3)$  succeeds with high
probability (whp), and if $\gl >1.275$ then it fails whp.
Mitzenmacher \cite{M} used the differential equation method
introduced by Wormald \cite{W} to claim that the threshold for the
PLA exists and it is the solution of certain equations, which are
somewhat complicated. That is, there is $\gl(k)$, $k\geq 3$, so
that the PLA applied to $F(n, \frac{\gl}{n^{k-1}}\,;k)$ succeeds
whp if $\gl< \gl(k)$, and fails whp if $\gl > \gl(k)$. It,
however, remains unclear whether it should be regarded as a
rigorous proof.

A more advanced algorithm called the unit clause  algorithm (UCA)
and its variations are analyzed \cite{CF,Dsur,A1,AS} to eventually
obtain the lower bound of 3.26. The UCA  first chooses a literal
uniformly at random and set it true. Then the negation of the
literal  is removed from the clauses containing it so that they
become  a clause of length one less. If there are clauses of
length $1$, or unit clauses, then the UCA chooses a  clause
uniformly at random among all unit clauses and set the literal in
the chosen clause true. The negation of the literal is removed
from the clauses containing it. Thus, it is possible that a
$0$-clause, i.e., a clause without any literal, can be created.
The UCA succeeds if no $0$-clause is created.

In a recent paper \cite{PCM1}, the author introduced the Poisson
cloning model $F_{PC} (n,p\,; k)$ for random $k$-SAT formulae,
which is essentially equivalent to the classical model $F (n,p\,
;k )$ when $p=\Theta (n^{1-k})$. That is,
\begin{thm} \label{equiv}
Let $k\geq 2$ and $ p= \Theta (n^{1-k})$.  Then there are
constants $c_{_1}$ and $c_{_2}$ such that, for any collection
${\cal F}$ of $k$-SAT formulae,
 $$ c_{_1}  \pr [
\fnpk \in \f] \leq    \pr [ F (n,p\,;k) \in \f]
 \leq
c_{_2}  (\pr [ \fnpk \in \f ]^{{\frac{1}{k}}} + e^{-n}),
 $$
 where
$$c_{_{1}}  = k^{1/2}e^{\frac{p}{n}{k \choose 2}{2n\choose k}
+ \frac{p^2}{2}{2n\choose k}}+o(1), ~~c_{_{2}}  =
e^{\frac{p(1-1/k)}{2n}{k \choose 2}{2n\choose
k}}\bb(\frac{k}{k-1}\bb) \bb((k-1)c_{_{1}}  \bb)^{1/k}+o(1).
$$
and $o(1)$ goes to $0$ as $n$ goes infinity.
  \old{
Moreover, if
$\pr [ \fnpk \in \f ] \geq n^{-b}$ for a constant  $b>0$ , then
$$\pr [ F (n,p\,;k) \in \f]
 \leq
c\pr [ \fnpk \in \f ],$$ for some constant $c$ depending on $b$. } 

\end{thm}

The  cut-off line algorithm (COLA) for  the new model is also
introduced in a general framework. Using the  COLA, one may
 generate an instance of the Poisson cloning model and
simultaneously carry an algorithm such as  the PLA. A version of
the  COLA applied to $F_{PC}(n,p\,;k)$ is analyzed to obtain the
following result for $F (n,p\, ; k)$: Let
$$ \gl (k)  := \min_{\rho>0} \frac{\rho
    }{(1-e^{-\rho})^{k-1}}, $$
    and $F_C (n,p)$ be the residual formula left after the PRA applied
$F(n,p)$ stops. The residual formula is called the core of
$F(n,p)$. The set of underlying variables of  $F_C(n,p)$ is
denoted by     $C (n,p\,;k)$. In other words, a variable is in $
 C(n,p\,;k)$  if and only if a clause of $F_C(n,p\,; k)$ contains
 it.

\begin{thm}
\label{plrthm} Let $\gl(n,p\, ;k)= p{2n-1 \choose k-1}$,
$k\geq 3$ and $\gs\gg n^{-1/2}$. 
Supercritical Phase:  If $\gl(n,p\, ;k)  <  \gl(k) -\gs $ is
uniformly bounded from below by $0$ and $i_{_0}(k)$ is the minimum
$i$ such that $2^k {i \choose k} \geq 2i/k$, then
$$ \Pr[ C (n,p\,;k) \not=\emptyset \, ] \leq 2 e^{-\Omega(\sigma^2 n)}
+ O(n^{-(1-2/k)i_{_0}(k)}).$$ \mn Supercritical Phase: If $\gl:=
\gl(n,p\, ;k) =\gl(k) +\gs $ is uniformly bounded from above,
then, for  the largest solution $\thl$ of the equation
$\th^{\frac{1}{k-1}} -1+e^{-\th\gl} =0$ and all $\ga $ in the
range
 $1\ll \ga  \ll \gs n^{1/2} $,
 $$  \pr [\, \, |\, |C(n,p\, ; k)|- \thl^{\frac{2}{k-1}} n | \geq
 \ga  ( n/\gs)^{1/2}\,  ] = e^{-\Omega (\ga ^2 ) }.$$
 In particular, the PRA succeeds with high probability if $\gl(n,p\, ;k)
 = \gl(k) -\gs$ with $\gs\gg n^{1/2}$, and it does not
  succeed with high probability if
$\gl(n,p\, ;k)
 = \gl(k) + \gs$ with  $\gs\gg n^{1/2}$.
\end{thm}
Most of  structural properties of the core can be found in
\cite{PCM1} too. The Poisson cloning model and the cut-off line
algorithm will be presented in detail in the next section .

\bn

For $k=2$, the PLA may not  succeed with nontrivial probability
even for $\gl_p :=p(2n-1) < 1$. For example,  there could be a
pair of clauses $(y\vee z) $ and $(\bar{y}\vee \bar{z})$ for two
variables $y$ and $z$ with non-trivial probability. Hence, we may
expect, at best, that if $\gl_p:= p(2n-1) <1$ then $F_C
(n,p):=F_C(n,p\,;2)$ consists of variables of type $(1,1)$ only.
Here and  in general, a variable $x$ is of type $(i,j)$ in a
formula if $x$ appears in $i$ clauses and $\bar{x}$ appears in $j$
clauses of the formula. The type of a literal $\bar{x}$ is
determined by the type of $x$. Taking similar approaches used to
analyze the structure of the core of the random digraph
\cite{PCD},  we will actually prove it and, in case that
 $\gl_p
>1$, we prove  that
 $F_C (n,p)$ has many   variables of type larger $(1,1)$ and
 the formula is
not satisfiable whp. All the proofs presented here do not depend
on \cite{PCD} though.

Other interesting properties for $F_C(n,p)$ are studied too.
Denoted by $C_{n,p}(i,j)$ is  the set of all variables of type
$(i,j)$ in $F_C(n,p)$ and $C_{n,p} = \cup_{(i,j)\geq (1,1)}
C_{n,p}(i,j)$ is the set of underlying variables of $F_C(n,p)$.
Due to the following lemma, the structure of the core $F_C(n,p)$
can be well understood provided  tight  upper and lower bounds for
$|C_{n,p}(i,j)|$'s are found, $(i,j) \geq (1,1)$.

\begin{thm}\label{config}
Suppose  two formulae 
have the same number of clauses on the same number of underlying
variables, and all underlying variables are of type at least
$(1,1)$. Then the two formulae are equally likely to be the core
of $F(n,p)$.
\end{thm}

\vspace{-3mm} \noindent The proof of the theorem is not difficult
and presented in Section \ref{sw2}.

For variables $x$ of type $(1,1)$ in $F_C(n.p)$, the conjunction
$(x\vee y) \wedge(\bar{x} \vee z)$ of  two clauses containing $x$
and $\bar{x}$ may be replaced by $(y\vee z)$. The replacement is
called a resolution of $x$. It is clear that the satisfiability is
not affected by a series of such resolutions. The formula obtained
after all possible resolutions of type $(1,1)$ variables is called
the kernel of $F(n,p)$ and denoted by $F_K (n,p)$. It is worth to
notice that  clauses in $F_K(n,p)$ may not consist of strictly
distinct literals.   Clearly, all variables of $F_K(n,p)$ are of
type larger than $(1,1)$, counting a loop $(x\vee x)$ twice in the
degree of $x$.

Let $D_{n,p} (i,j)= \sup_{(i',j')\geq (i,j)} C_{n,p} (i',j')$.
Then $D_{n,p} (1,1) = C_{n,p} $ and $K_{n,p}:=D_{n,p} (2,1) \cup
D_{n,p} (1,2)$ is the set of underlying variables of $F_K (n,p)$.
When $C_{n,p} (i',j')$ are all small for $(i',j') \geq (i,j)$, it
sometimes more useful and/or easier to bound the size of $D_{n,p}
(i,j)$ rather than individual $C_{n,p} (i',j')$.
 We also denote  $M_{n,p} (i,j)$ to be the sum of
degrees of all variables in $D_{n,p} (i,j)$ and their negations,
where the degree $d(y)$ of a literal $y$  is the number of clauses
 containing it. Clearly,
$$ M_{n,p} (i,j) = \sum_{(i',j') \geq (i,j)} (i'+j') |C_{n,p}
(i',j')|. $$ Notice that  the numbers of clauses in $F_C (n,p)$
and $ F_K(n,p)$ are $\frac{1}{2}M_{n,p} (1,1)$ and
$\frac{1}{2}(M_{n,p} (1,2)+M_{n,p} (2,1)-M_{n,p} (2,2))$,
respectively. Finally, we set
$$P_\ell(\mu) = \pr[ \poi (\mu)= \ell] = e^{-\mu}
\frac{\mu^\ell}{\ell!}, ~~{\rm and} ~~Q_\ell(\mu ) = \pr[ \poi
(\mu)\geq  \ell] = e^{-\mu} \sum_{\ell'\geq \ell}
\frac{\mu^{\ell'}}{\ell'!}. $$

In statements in theorems, lemmas and corollaries of this paper,
we use the following convention.

\mn {\bf Convention:} When we say that a statement is true for all
$\ga$ in the range $a\ll \ga \ll b$, it actually means that there
is (small) constant $\eps>0$ so that the statement is true for
$\ga$ in the range $a/\eps \leq \ga \leq \eps b$.

\old{
 A variables $x_i$  is called {\em type $(a,b)$} in a formula
if the degree of $x_i$ is $a$ and the degree of $\bar{x}_i$ is
$b$. If a variable is of type $(1,1)$, the variable and the
corresponding literals are called {\em neutral}, as these literals
have no effect to the formula after the resolution of $(y\vee z)
\wedge (\bar{z}, y')$ into $(y\vee y')$. A literal that is neither
pure nor neutral is called {\em complex}.
 } 

\begin{thm}\label{seccor}   Suppose  $p(2n-1) =1+\gs $ is uniformly bounded from
above with $\gs \gg n^{-1/3}$. Let $\gl=1+\gs$, $\gD>0$ and $1\ll
\ga\ll (\thl n)^{1/2} $. Then, for  fixed  $(i,j)$ and  $P_i = P_i
(\thl \gl)$ and  $Q_i = Q_i (\thl \gl)$,
$$\pr\bb[\,\, \bb| |C_{n,p} (i,j)| -  P_i P_j  n
\bb| \geq   \ga \thl^{i+j-3} (\thl^3 n)^{1/2}  \bb] \leq
e^{-\Omega(\thl^{i+j-3} \ga^2)}+ e^{-\Omega(\ga^2)},$$ and,
assuming $i\geq j$,
$$\pr\bb[\,\, \bb| |D_{n,p} (i,j)\cup D_{n,p} (j,i)|
- \bb( 2 Q_i Q_j -  Q_i Q_i \bb) n\bb| \geq \ga \thl^{i+j-3}
(\thl^3 n)^{1/2} \bb] \leq e^{-\Omega(\thl^{i+j-3} \ga^2)}+
e^{-\Omega(\ga^2)},$$ and
$$
\pr\bb[\,\, \bb|M_{n,p} (i,j) -  \thl\gl \bb(Q_{i-1}  Q_j + Q_{i}
Q_{j-1} \bb) n \bb| \geq \ga \thl^{i+j-3} (\thl^3 n)^{1/2} \bb]
\leq e^{-\Omega(\thl^{i+j-3} \ga^2)}+ e^{-\Omega(\ga^2)}.$$
Moreover, $$ \pr [ |D_{n,p} (i,j) | \geq \ell ] \leq O
\bb(\frac{((1+\frac{\ga}{(\thl n)^{1/2}})\thl
\gl)^{(i+j)\ell/2}}{(\ell!)^{1/2}}\bb) + e^{-\Omega (\ga^2)}.
$$
   \old{
Moreover, for $(i,j)\geq (2,2)$ and, for  any $\gd$ in the range
$\thl^{i+j-3}\leq \gd \ll 1$,
$$\pr \bb[ |D_{n,p} (i,j)| \geq \sfrac{(\thl\gl)^{i+j}n}{i!j!}
+ \gd \thl^3 n \bb] \leq e^{-\Omega(\gd \thl^3 n )},
$$
and
$$\pr \bb[ M_{n,p} (i,j) \geq \sfrac{(\thl\gl)^{i+j}n}{(i-1)!j!}
+ \sfrac{(\thl\gl)^{i+j}n}{i!(j-1)!} + \gd \thl^3 n \bb] \leq
e^{-\Omega(\gd \thl^3 n)}.
$$} 
\end{thm}

\mn

A stronger theorem (Theorem \ref{main}, see also Main Lemma in
Section \ref{sw}) is to be first proved and Theorem \ref{seccor}
will follow as a corollary. Bounds for the sizes of the core and
the kernel may be obtained from Theorem \ref{seccor}. Estimations
for $|F_C(n,p)|, |F_K(n,p)|$ are possible too, where, in general,
$|F|$ is the number of clauses in the formula $F$.

\begin{cor}\label{corekernel} For the core $F_C(n,p)$ of $F(n,p)$ and
 the  set
$C_{n,p}$ of underlying variables of  the core,
$$ \pr \bb[\, \,\bb| |C_{n,p}| - \thl^2 n
\bb| \geq \ga (\thl n)^{1/2}\bb]
 \leq e^{-\Omega(\ga^2)}, $$
 and
$$ \pr \bb[\, \,\bb| |F_C(n,p) | -\thl^2\gl  n
\bb| \geq \ga (\thl n)^{1/2}\bb]
 \leq e^{-\Omega(\ga^2)}. $$
For the kernel $F_K(n,p)$ of $F(n,p)$ and the set $K_{n,p}$ of
underlying variables of the kernel,
$$ \pr \bb[\, \,\bb| |K_{n,p}| -\thl^2
(1-\gl^2 e^{-2\thl\gl} )n \bb| \geq \ga (\thl^3 n)^{1/2}\bb]
 \leq e^{-\Omega(\ga^2)}, $$
$$
\pr \bb[\, \,\bb| |F_K(n,p)| -\thl^2 \gl (1-\gl e^{-2\thl\gl} )n
\bb| \geq \ga (\thl^3 n)^{1/2}\bb]
 \leq e^{-\Omega(\ga^2)}. $$
   \old{
$$ \pr \bb[ M_{n,p} (2,2) \geq (\thl \gl)^4 n + \gd \thl^3 n
 \bb]
 \leq e^{-\Omega(\gd \thl^3 n)}.
$$} 
\end{cor}

In brief, we may also have
\begin{cor}
\label{swthm}  Let  $\gl=1+\gs$ with  $n^{-1/3} \ll \gs < 1$.
Then, with high probability, the pure literal algorithm  applied
to $F(n,\frac{\gl}{2n-1})$ stops leaving $\Theta(\gs^2 n)$ type
$(1,1)$ variables, $\Theta(\gs^3 n)$ type $(2,1)$ or $(1,2)$
variables, and $O(\gs^4 n)$ clauses containing  other type
variables. Moreover, once $C_{n,p} (i,j)$, $(i,j)\geq (1,1)$, are
given, the residual formula is the uniform random formula
conditioned on $C_{n,p} (i,j)$.
\end{cor}

The analysis of the structure of the core yields almost optimal
bounds  for the probability of  satisfiability, improving bounds
of Bollob\'as et. al. \cite{BBCKW}.

\begin{thm}
\label{subker} If $\gl_p =1 -\gs $ is uniformly bounded from below
by $0$ with $n^{-1/3}\ll \gs \ll 1 $, then, with probability
$1-\frac{15+o(1)}{16\gs^3 n}$, all the variables in $C (n,p)$
 are  of type $(1,1)$. That is,
  \begin{equation}\label{supersize} \pr \bb[\, \,
  K(n,p)= \emptyset
\, \, \bb] = 1-\frac{15+o(1)}{16\gs^3 n}.
 \end{equation}
 In particular, $\pr [K(n,p)= \emptyset
\, ] = 1-O((\gs^3 n)^{-1})$ for all $\gs$ in the range $n^{-1/3}
\ll \gs <1$.
 We also have
 $$ \pr [\,  F (n,p) ~\mbox{is satisfiable}\, \, ]
 = 1-  \frac{1+o(1)}{16\gs^3 n}. $$
\end{thm}

\begin{thm}
\label{sub2sat} If $\gl_p =1 +\gs $ is uniformly bounded from
above, then $F(n,p)$ is unsatisfiable  with probability
$1-e^{-\Theta(\gs^3 n)}$, i.e.,
 $$ \pr [\,  F (n,p) ~\mbox{is satisfiable}\, \, ]
 = e^{-\Theta(\gs^{3} n)}. $$
\end{thm}

In the next section, we present the Poisson cloning model and the
cut-off line algorithm together with an useful large deviation
inequality called generalized Chernoff bound. Then, Theorem
\ref{seccor} and Corollaries \ref{corekernel} and  \ref{swthm}
will be proven in Section \ref{sw}. Section \ref{sw2} is for the
proofs of Theorems \ref{config} \ref{subker} and \ref{sub2sat}.

\include{dstart}
\newcommand{\nopf}[1]{}

\dvi

\section{Poisson Cloning Model and Cut-Off Line Algorithm} 
\label{model}

\newcommand{\bin}{{\rm Bin}}

{\bf Poisson Cloning Model:} The Poisson cloning model is
partially motivated by the fact that the degree $d(y)$ of a
literal $y$ in $F(n,p)$ is the binomial distribution $\bin (2n-1,
p)$, which is close to $\poi(p(2n-1))$ when $p=\Theta (n^{-1})$.
Here the degree $d(y)$ of $y$ is the number of clauses in $F(n,p)$
containing $y$ and
$$ \pr[ \bin (2n-1, p)=\ell] = {2n-1 \choose \ell} p^\ell
(1-p)^{2n-1-\ell}, ~~~\pr[ \poi(\gl)=\ell] =e^{-\gl}
\frac{\gl^\ell}{\ell!} .
$$
Though the degrees $d(y)$'s are not exactly independent, they are
expected to behave like i.i.d random variables. Thus, it has been
desirable to introduce a new model for the random $k$-SAT formulae
in which the degrees are i.i.d Poisson random variables. Inspired
by the configuration model for random regular graphs, see e.g.
\cite{BBK}, \cite{BC}, \cite{Bo}, and \cite{Wo}, the author have
introduced the Poisson cloning model with the desired properties
and show that the new model is not much different from the
classical model in the sense of Theorem \ref{equiv}.

To analyze various properties of random graphs and random SAT
formulae such as cores and giant components, the cut-off line
algorithm is introduced too. In this section, we present the
Poisson cloning model and the cut-off line algorithm, and related
lemmas as well as a large deviation inequality called generalized
Chernoff bound.

For  a new random $2$-SAT model $\fnp$, we take i.i.d Poisson
$\gl:=p(2n-1)$ random variables $d_y$ for each $y$ in the set $Y$
of all literals, and then take $d_y$ copies of each $y$. The
copies of a literal $y$ are called {\em clones of $y$}, or simply
{\em $y$-clones}. Since the sum of Poisson random variables is
also Poisson, the total number $N_\gl := \sum_{y\in Y} d_y$ of
clones is a Poisson $2\gl n$ random variable. It is sometimes
convenient to take a reverse, but equivalent, construction. We
first take a Poisson $2 \gl n$ random variables $N_\gl$ and then
take $N_\gl$ unlabelled clones. Each clone is independently
labelled as $y$-clone uniformly at random, in the sense that $y$
is chosen uniformly at random from $Y$. It is well-known that the
numbers $d_y$ of $y$-clones  are i.i.d Poisson $\gl$ random
variables.

 If $N_{\gl}$ is even, the formula $\fnp $ is to be
defined by generating a (uniform) random perfect matching on those
$N_\gl$ clones and contracting  clones of a literal $y$ into $y$.
That is, an edge consisting of a $y$-clone and a $z$-clone in the
perfect matching yields the clause $(y \vee z)$ in $\fnp$ with
multiplicity. If $y=z$, it produces, a loop $(y\vee y)$, which
contributes $2$ in the degree of $y$. It turns out that there are
many ways  to generate the random perfect matchings and we may
choose one that makes given problems easier to analyze.  Some
specific ways will be discussed when the cut-off line algorithm is
introduced.

 If $N_\gl$ is odd, we arbitrarily choose a clone, say $y$-clone. This
  clone induces
 a $1$-clause, called a defected
clause, consisting of $y$. The defected clause contribute only 1
to the degree of the corresponding literal. The same procedure
taken for the case of even $N_\gl$ are to be carried for the rest
of clones. Strictly speaking $\fnp$ varies depending on how to
construct the defected clause. However, for any collection
$\mathcal{F}$ of $2$-SAT formulae, the probability that $\fnp$ is
in $\mathcal{F}$ does not depend on how the defected clause is
chosen (for odd $N_\gl$),  since $\fnp \not\in \mathcal{F}$
whenever there is a non-standard clause in $\fnp$. Thus it is
normally unnecessary to describe $\fnp$ for odd $N_{\gl}$. For
$k\geq 3$, the Poisson cloning model $\fnpk$ for random $k$-SAT
problems may   be similarly defined.

Theorem \ref{equiv} has been proved  using somewhat
straightforward computations for $\pr[F(n,p \,; k)=F]$ and $ \pr[
\fnpk = F]$.

\bn

\bn {\bf Cut-Off Line Algorithm (COLA):}  To generate a uniform
random perfect matching on $N_\gl$ clones, we may keep matching
two unmatched clones uniformly at random.  Another way is to
choose the first clone as we like and match it to a clone chosen
uniformly at random among all other unmatched clones. Clearly,
there are many ways to choose the first clone. This is a big
advantage since we may select a way that makes the given problem
easier to analyze. In general, a sequence of choice functions will
tell how to choose the first clone at each step.  A choice
function may be deterministic or random. If $N_\gl$ is even, this
would yield a uniform perfect matching regardless what the choice
functions are. If only one clone, say of $y$,  remains unmatched,
we just add the defected clause consisting of $y$.

It is  useful to introduce a more specific way to choose the
second clone  uniformly at  random. The way presented here will be
useful to analyze some algorithms like the  PRA. First, we
independently assign, to each clone, a uniform random real number
between $0$ and $\gl$. For the sake of convenience, we say that a
clone is the largest, smallest, etc. if so is its assigned number.
Each choice function is to choose an unmatched clone without
changing the (joint) distribution of the numbers assigned to all
other unmatched clones. A choice function satisfying this
condition is called {\em oblivious}. For instance, a choice
function is oblivious if it
 chooses a clone of a pure literal. If a
choice function chooses a largest  $v$-clone, it is not oblivious,
as it changes the distribution of the numbers assigned to other
unmatched $v$-clones.

   Once
an unmatched clone is chosen by an oblivious choice function, the
largest clone among all other unmatched clones are  to be matched
to the chosen clone. This may be further implemented using the
Poisson $\gl$-cell: First, map a $y_j$-clone with assigned number
$r$ to  the point $(r, j)$ in the two dimensional plane.  One may
think that there are $2n$ horizontal line segments  in
$\mathbf{R}^2 $ from $(0,j)$ to $(\gl, j)$, $j=1, ..., 2n$ and, on
each line segment, there are i.i.d. uniform $d_{y_{_j}}$ points
that tell the assigned numbers for $d_{y_{_j}}$ clones of
$y_{_j}$. This rectangular configuration is called a {\em Poisson
$\gl$-cell}. Each line segment of the Poisson $\gl$-cell with the
points  is an independent Poisson arrival process with density
$1$, up to time $\gl$.

The cut-off line algorithm (COLA) can be described as follows.
Initially, the cut-off line is the vertical line in $\mathbf{R}^2
$ containing the point $(\gl, 0)$. At the first step, once the
oblivious choice function chooses a clone, we move the cut-off
line to the left until a clone is on the line. The clone is
clearly the largest unmatched clone, excluding the chosen clone.
The new cut-off value, denoted by $\gL_1$, is the assigned number
to the clone. The new cut-off line is, of course, the vertical
line containing $(\gL_1, 0)$. Keep repeating this procedure, one
may obtain the $i^{\rm th}$ cut-off value $\gL_i$ and the
corresponding cut-off line. It is crucial to note that, provided
all choice functions are oblivious, once $\gL_{i}$ is given then
all numbers assigned to unmatched clones are i.i.d uniform random
numbers between $0$ to $\gL_i$.

For $\th$ in the range $0\leq \th \leq 1$, let $\gL(\th)$ be the
cut-off value when $(1-\th^{2}) \gl n $ or more clones are matched
for the first time. Conversely, let $N(\th)$ be the number of
matched clones until the cut-off line reaches  $\th \gl$. Two
versions of the cut-off line lemma have been proven in
\cite{PCM1}.

\begin{lem} (Cut-off Line Lemma)
\label{cf2} Let  $\gl >0$ be fixed.
 Then, for $\tho <1 $ uniformly bounded below from $0$
 and $0<\dd\leq n$,
$$ \pr \bb[  \max_{\th:\tho \leq \th \leq 1 } |\gL (\th) -\th \gl|  \geq
\sfrac{\dd}{n}
 \bb]
  \leq  2e^{-\Omega(\min\{\dd, \frac{\dd^2}{(1-\tho)n}\})},$$
and $$ \pr \bb[  \max_{\th:\tho \leq \th \leq 1 } |N(\th)
-2(1-\th^\kk) \gl n|
 \geq \dd  \bb]
  \leq  2e^{-\Omega(\min\{\dd, \frac{\dd^2}{(1-\tho)n}\})}.
    $$
\end{lem}

For the Poisson $\gl$-cell conditioned on $N_\gl=N$, a similar
lemma may be obtained.

\begin{lem} (Cut-off Line Lemma for $N$ clones)
\label{cf3} Let $k\geq 2$, $\gl >0$ be fixed.
 Then, for the Poisson $\gl$-cell conditioned  on $N_\gl=N$, and for
 $\tho <1 $ uniformly bounded below from $0$
 and $0<\dd\leq N$,
$$ \pr \bb[  \max_{\th:\tho \leq \th \leq 1 } |\gL (\th) -\th \gl|  \geq
\sfrac{\dd}{N}
 \bb]
  \leq  2e^{-\Omega(\min\{\dd, \frac{\dd^2}{(1-\tho)N}\}},$$
and $$ \pr \bb[  \max_{\th:\tho \leq \th \leq 1 } |N(\th)
-(1-\th^\kk) N|
 \geq \dd  \bb]
  \leq  2e^{-\Omega(\min\{\dd, \frac{\dd^2}{(1-\tho)N}\}}.
    $$
\end{lem}

For the proof of the cut-off line lemma, a  large deviation
inequality, called generalized Chernoff bound, has been used.
Here, we present a version of it  that is useful for our analysis.
A proof can be found in \cite{PCM1}.

\newcommand{\vp}{\xi}

\begin{lem} \label{uni}(Generalized Chernoff bound)
Let $ X_1, ..., X_m $ be a sequence of   random variables. Suppose
  \beq{mean} E[X_i| X_{1}, ..., X_{i-1} ] \leq
\mu_i , \eeq and there are  $a_{_i}$, $b_{_i}$ and $\xi_{_0}$ so
that
 \beq{second}  E[ (X_i-\mu_i)^2|X_{1}, ..., X_{i-1} ]
\leq a_i, \eeq
 and
  \beq{thirdco}
E[ (X_i - \mu_i)^3 e^{\vp (X_i-\mu_i) }| X_{1}, ..., X_{i-1} ]
\leq b_i ~~~\mbox{for all $0\leq \vp \leq \vp_{_0} $}.\eeq
  If   $ \gd\vp_{_0} \sum_{i=1}^m b_i \leq  \sum_{i=1}^m a_{_i} $
for some $0<\gd\leq 1$, then
 $$ \pr \bb[ \sum_{i=1}^m X_i   \geq  \sum_{i=1}^m \mu_i+\gD  \bb] \leq
e^{-\frac{1}{3}\min \{\gd  \vp_{_0}\gD, \, \,  \gD^2/ \sum_{i=1}^m
a_i
                        \})},
$$ for all $\gD>0$.
Furthermore, if  $ X_1, ..., X_m $ are independent and  satisfy
\raf{second} for $\mu_i= E[X_i]$ and
  \beq{thirdco2}
\bb| E[ (X_i - E[X_i] )^3 e^{\vp (X_i-E[X_i]) } ] \bb| \leq b_i
~~~\mbox{for all $\vp$ in the range $|\vp | \leq \vp_{_0} $},\eeq
then  $ \gd\vp_{_0} \sum_{i=1}^m  b_i \leq  \sum_{i}^m a_{_i} $
for $0<\gd\leq 1$ implies that
$$ \pr \bb[ \,\, \bb|\sum_{i=1}^m X_i   -  \sum_{i=1}^m E[X_i]\bb| \geq \gD
  \bb] \leq e^{-\frac{1}{3}\min \{\gd  \vp_{_0}\gD , \, \,
  \frac{\gD^2}{
\sum_{i=1}^m a_i}
                        \})},
$$ for all $\gD>0$.
\end{lem}

We conclude this section by presenting a  corollary that can be
applied to random walks with negative drift.

\begin{cor}\label{nd}
Suppose \raf{mean}-\raf{thirdco} hold with $\mu_i = -h$ for a
constant  $\gb >0$.
 If   $ \gd\vp_{_0} \sum b_i \leq  \sum a_{_i} $
for some $0<\gd\leq 1$, then
 $$ \pr \bb[ \sum_{i=1}^m X_i   \geq  \gD  \bb] \leq
e^{-\Omega(\min \{\gd  \vp_{_0}(\gD+hm)  , \, \, (\gD+hm)^2 /
\sum_{i=1}^m a_i\})}.
$$
\end{cor}

\include{dstart}

\dvi

\newcommand{\bv}{\bar{v}}
\newcommand{\bx}{\bar{x}}

\newcommand{\smp}{e^{-\Omega(\ga^2)}}
\newcommand{\smpr}{e^{-\Omega_\gamma (\ga^2)}}

\newcommand{\hh}{H_{PC}}
\newcommand{\hl}{H_{PC}(\gl)}

\include{dstart}

\topmargin 0pt \advance \topmargin by -\headheight \advance
\topmargin by -\headsep \textheight 8.9in \oddsidemargin 0pt
\evensidemargin \oddsidemargin \marginparwidth 0.5in \textwidth
6.5in

\dvi
\section{Pure literal algorithm   for the random  $2$-SAT problem}
\label{sw}

As mentioned in the previous section, the COLA is useful to
realize some algorithms like the PLA. The following specific COLA
is used to analyze the structure of the core of $F_{PC} (n,p)$.

 \mn {\bf COLA (for core)}: Construct a Poisson $\gl$-cell.
 If a variable is  of types $(0,i)$ or $(i,0)$, put
 all clones of it and its negation into  a
stack in an arbitrary order. This does not mean that the clones
are removed from the $\gl$-cell.

\mn (a)   If the stack is empty, go to (b). If the stack  is
nonempty, choose the first clone in the stack and move the cut-off
line to the left  until the largest  unmatched clone, excluding
the chosen clone, is found. (The stack naturally defines  choice
functions.) Then, match the largest unmatched clone to the chosen
clone. Remove all matched clones from the stack  and from the
cell. If there are new variables of type $(0,i)$ or  $(i,0)$, then
put all clones of them and their negations in the stack. Repeat
(a).

\mn (b) Choose a clone uniformly at random from all unmatched
clones and put it in the stack. Then, go to (a).

\mn The steps carried by the instruction described in (b) are
called {\em free steps} as it is free to choose any clone. We will
call unmatched  clones of pure literal {\em light} and the other
unmatched  clones {\em heavy}. A literal is called heavy if it not
pure.

According to the cut-off line lemma, one may expect that there are
$2\th^2 \gl n$ unmatched clones (when the cut-off line is) at $\th
\gl$. The number of heavy clones at $\th\gl$ is expected   to be
close to  $2(1-e^{-\th\gl})\th\gl n=\Theta(\th^2 n)$. (See
\raf{v4hh} below.) Thus, the number of light clones seems to be
close to
 $$2\th^2\gl n - 2(1-e^{-\th\gl})n = 2\th\gl n
 (\th-1+e^{-\th\gl}), $$
 which is $\Theta(\th^3 n)$ provided $\th \gg |\gl-1|$.
 If $\th$ is small, however, this observation
  would give us no information. This is  due to the fact that
  the standard deviation for the number
of heavy clones is $\th n^{1/2}$ so that, for $\th^3 n \ll \th
n^{1/2}$, or $\th \ll n^{-1/4}$,  it is unclear wether the number
of light clones is positive or not.

A more careful analysis starts from the observation that, when
$\th$ is small, most of heavy variables are of type $(1,1)$ and
that the two clauses
 containing such a variable and its negation may
 be resolved to one clause.
In other words, the two clause $(x\vee y)$ and  $(\bar{x}\vee z)$
may be replaced by
 $(y\vee z)$, which is called a resolution. After a
series of such resolutions, all variables of type $(1,1)$ may
disappear.

To take an advantage of this fact, we will introduce many phases.
Let $\frac{1-\thl}{10}\leq \gb\leq \frac{1-\thl}{2}$. The first
phase starts at the beginning of the whole process. For $j\geq 1$,
the $j^{\rm th}$ phase ends and the $(j+1)^{\rm th}$ phase begins
when the cut-off line reaches $(1-\gb)^{j} \gl$. At the beginning
of each phase, all variables of type $(1,1)$ and their unmatched
clones are called passive. All other unmatched clones are called
active. These terms do not change until the beginning of the next
phase. So, variables that become type (1,1) only after the current
phase starts remain active until the end of the phase. Once  a
clone becomes pure, it plays the same role regardless of being
passive or active. The procedure (b) of COLA  also need to be
replaced by

\mn (b)* Choose a clone uniformly at random from all unmatched
active clones and put it in the stack. If there is no active
clone, stop. Otherwise,  go to (a).

As a stack is used, if  one of the two unmatched clones of  a
passive variable and its negation were matched in a step then  the
choice function in the next step must  choose the other clone.
Thus, the situation is exactly the same except  the number of
passive variables decreases by $1$. This means that the COLA
 applied
without passive clones  is essentially the same as the original
algorithm.
 In this sense, we may say that two
active clones are matched  if so are they after the resolutions of
matched passive clones.
 Here the resolution has the natural meaning: Two
edges $\{z_1, z_2\}$, $\{z_3, z_4\}$ with clones $z_2, z_3$ of a
passive variable and its negation is reduced to the one edge
$\{z_1, z_4\}$. Conversely, an active clone may be regarded as
unmatched if it is not matched or it is not matched after the
resolutions.

Let $\gL_C$ be the cut-off value when no light clone remains for
the first time in the COLA applied to the Poisson $\gl$-cell. The
main lemma shows that $\gl_C$ is highly concentrated near $\thl
\gl$, as expected,  with standard deviation $(\thl n)^{-1/2}$.
Once $\gL_C$ is determined, the unmatched clones  form the Poisson
 $\gL_C$-cell without pure literals.

\begin{lem}\label{upper} (Main Lemma)
Let $\gl =1+\gs$ with $\gs \gg n^{-1/3}$. Then, for all $\ga$ with
$1\ll \ga \ll (\thl^3 n)^{1/2}$,
$$ \pr  [|\gL_C  - \thl \gl| \geq \ga (\thl n)^{-1/2} ]
=e^{-\Omega(\ga^2)}. $$
\end{lem}

\bn

 For the proof, we first  estimate the number of
active clones at the beginning of each phase. Let $N_j$ be  the
number  of active clones at the beginning of the $j^{\rm th}$
phase and let $M_j$ be the number  of matched active clones during
the entire  $j^{\rm th}$ phase. Then, the cut-off line lemma for
$N_j$ clones, or Theorem \ref{cf3},  gives \beq{v3mm} \pr [ \, \,
| M_j - (1-(1-\gb)^2) N_j |\geq \gD| N_j ] \leq 2
e^{-\Omega(\min\{\gD, \frac{\gD^2}{N_j}\}) }. \eeq
  Notice that the number $N_{j+1}$ of
active clones at the beginning of the next phase is $N_j - M_j -
2B_j$,
where $B_j$ is the number of variables of type  $(1,1)$ at
$(1-\gb)^j \gl$ that were  of type larger than $(1,1)$ at
$(1-\gb)^{j-1} \gl$. (Recall that an active clone is regarded as
unmatched if it is not matched or it is not matched after
resolutions.)  For a literal $y$ and $0\leq \th < \th' \leq 1$,
denoted by  $d_y (\th, \th')$ is the number of $y$-clones larger
than or equal to $\th\gl$ and smaller than $\th'\gl$, and $d_y
(\th) = d_y (0, \th)$.
Then, for $\th_j = (1-\gb)^{j-1}$,
$$ B_j =\sum_{x\in X}  1(d_x (\th_{j+1})=d_{\bar{x}}
(\th_{j+1})=1) 1(d_x(\th_{j+1}, \th_j) + d_{\bar{x}}
(\th_{j+1},\th_j)\geq 1).
$$
Observe that $(d_x (\th_{j+1}),d_{\bar{x}} (\th_{j+1}),
d_x(\th_{j+1}, \th_j), d_{\bar{x}} (\th_{j+1},\th_j))$, $x \in X$,
are i.i.d $4$-tuples of independent Poisson random variables with
means $\th_{j+1} \gl$, $\th_{j+1} \gl$, $(\th_j-\th_{j+1}) \gl$,
$(\th_j- \th_{j+1}) \gl$, respectively. Applying the generalized
Chernoff bound, we have
  \beq{v3bb}\pr [
\, |B_j - (\th_{j+1}\gl)^2 e^{-2\th_{j+1} \gl} (1- e^{-2\gb \th_j
\gl})n| \geq  \gD ]  \leq 2e^{-\Omega(\min \{
\gD,\frac{\gD^{2^{\phantom{.}}}}{\th_j^3 n}\})}. \eeq Therefore,
$N_{j+1}$ is expected to be close to $(1-\gb)^2 N_j -
2(\th_{j+1}\gl)^2 e^{-2\th_{j+1} \gl} (1- e^{-2\gb \th_j \gl})n$.
Applying this inductively, we expect that $N_{j}$ is close to $
  2 \th_j^2\gl  ( 1-\gl e^{-2\th_j \gl })n,
$.

Let \bean H_j &=& \sum_{x\in X} (d_x(\th_j) + d_{\bar{x}} (\th_j))
1\bb( (d_x(\th_j), d_{\bar{x}} (\th_j))> (1,1)\bb)\\ &=&
\sum_{x\in X} (d_x(\th_j) + d_{\bar{x}} (\th_j)) 1\bb(
(d_x(\th_j), d_{\bar{x}} (\th_j))\geq (1,1)\bb) -2 \sum_{x\in X}
1\bb( d_x(\th_j) = d_{\bar{x}} (\th_j)=1\bb). \eean Then $H_j$ is
the number of active heavy clones at the beginning of the $j^{\rm
th}$ phase unless there is a free step before $\th_j \gl$.
Generally, $H_j$ is an upper bound for the number of heavy clones
and  $L_j:=N_j -H_j$ is a lower bound for the number of light
clones. The bounds may be strict only when there is a free step
before the cut-off line reaches $\th_j \gl$.

 As  $(d_x(\th_j), d_{\bar{x}} (\th_j))$ are i.i.d
pairs of independent Poisson random variables with mean $\th_j
\gl$, the generalized Chernoff bound gives
  \beq{v4hh} \pr \bb[ \, \,
\bb|H_j -2\th_j \gl (1-e^{-\th_j \gl}-\th_j\gl
 e^{-2\th_j\gl})n\bb| \geq \gD \, \bb]
 \leq 2e^{-\Omega(\min \{
\gD,\frac{\gD^{2^{\phantom{.}}}}{\th_j^3 n}\})}. \eeq

 Suppose $\gl = 1+ \gs $ with $\gs \gg n^{-1/3}$ and
$1\ll \ga \ll (\thl^3 n)^{1/2}$. We take $\frac{1-\thl}{10} \leq
\gb \leq \frac{1-\thl}{2}$ so that $ (1-\gb)^{a-1} = \thl + \ga
(\thl n)^{-1/2}$ for an integer $a$. Let
 $$\gD_j= 0.01\ga ( \th_j^3
n)^{1/2}  \sum_{i=1}^{j} (1-\gb)^{\frac{2j-i-a}{4}}. $$ Then,
since  $(1-\gb)^{\frac{j}{4} } \th_j^{-3/2} =(1-\gb)^{3/2}
(1-\gb)^{-5j/4}$ increase as $j$
increases, 
and $\th_a= \thl + \ga (\thl n)^{-1/2} = (1+o(1))\thl$, we have
\beq{v4aa}
  \ga (1-\gb)^{\frac{j-a}{4}} (\th_j^3 n)^{-1/2}   \leq \ga (\th_a^3
  n)^{-1/2} \ll 1 . \eeq
and
 $$ 
    \gD_j= 0.01\ga ( \th_j^3
n)^{1/2}  \sum_{i=1}^{j} (1-\gb)^{\frac{2j-i-a}{4}} = 0.01\ga (
\th_j^3 n)^{1/2}   (1-\gb)^{\frac{j-a}{4}} \sum_{i=1}^{j}
(1-\gb)^{\frac{j-i}{4}} \ll \th_j^3 n
$$
      for all $j=1,..., a$.
\begin{lem} \label{v3sup}  For all $\ell=1,...,a$, we have
$$ \pr \bb[ \, \exists ~j=1,...,\ell ~~s.t.~~
 |N_{j} - 2 \th_j^2\gl  ( 1-\gl e^{-2\th_j \gl })n| >  \gD_{j}   \bb] \leq
   e^{-\Omega(\ga^2(1-\gb)^{\frac{\ell-a}{2}})},
$$
and
$$ \pr \bb[ \, \exists ~j=1,...,\ell ~~s.t.~~
 |L_{j} - 2 \th_j\gl  ( \th_j- 1+ e^{-\th_j \gl })n| >  2\gD_{j}   \bb] \leq
   e^{-\Omega(\ga^2(1-\gb)^{\frac{\ell-a}{2}})}.$$

\end{lem}

\noindent \pf  Let $n_j = 2 \th_j^2\gl  ( 1-\gl e^{-2\th_j \gl
})n$, $b_j= (\th_{j+1}\gl)^2 e^{-2\th_{j+1} \gl} (1- e^{-2\gb
\th_j \gl})n$ and $\ga_j = (1-\gb)^{\frac{j-a}{4}} \ga$. Then
 $\ga_j \ll (\th_j^3 n)^{1/2}$ by \raf{v4aa} and
$$ n_{j+1} = (1-\gb)^2 n_j - 2b_j. $$
  Since
  \bean  \pr \bb[ \, \exists ~j=1,...,\ell+1 ~s.t.~
 |N_{j} -  n_j | >  \gD_{j}
 \bb]\! \! \! \! & =&  \! \! \! \pr \bb[ \, \exists ~j=1,...,\ell ~~s.t.~~
 |N_{j} -   n_j | >  \gD_{j}
 \bb] \\
 & & \! \!  \! \! \!
 + \pr \bb[
 |N_{j} -  n_j | \leq   \gD_{j}  ~\forall j\leq \ell,~
 |N_{\ell+1} - n_{_{\ell+1}}| >  \gD_{\ell+1}
 \bb]
\eean and
$$ \pr \bb[
 |N_{j} -  n_j | \leq   \gD_{j}  ~\forall j\leq \ell,~
 |N_{\ell+1} - n_{_{\ell+1}}| >  \gD_{\ell+1}
 \bb] \leq \pr \bb[
 |N_{\ell} -  n_{_{\ell}} | \leq   \gD_{\ell},~
 |N_{\ell+1} - n_{_{\ell+1}}| >  \gD_{\ell+1}
 \bb],
$$
it is enough by $\sum_{j=1}^{\ell+1} e^{-\Omega(\ga_{j}^2)}
 = e^{-\Omega(\ga^2_{{\ell+1}})}$ to show that
$$
P_{\ell+1} := \pr \bb[
 |N_{\ell} -  n_{_{\ell}} | \leq   \gD_{\ell},~~
 |N_{\ell+1} - n_{_{\ell+1}}| >  \gD_{\ell+1}
 \bb] \leq e^{-\Omega(\ga_{\ell+1}^2)}. $$

Notice that $N_{\ell+1} = N_\ell -M_\ell - 2B_\ell$,
$n_{_{\ell+1}} = (1-\gb)^2 n_{_\ell} - 2b_\ell $ and
 \bean |N_{\ell+1} - n_{_{\ell+1}} | &\leq & |N_\ell-M_\ell
  -(1-\gb)^2 N_\ell| + (1-\gb)^2 |N_\ell - n_{_\ell}|
+| (1-\gb)^2 n_{_\ell}  - n_{_{\ell+1}}-2B_\ell|\\
&=& |M_\ell - (1-(1-\gb)^2) N_\ell|+ (1-\gb)^2 |N_\ell -n_{_\ell}|
+2|B_\ell - b_{\ell}| .\eean
 As
 $ \gD_{\ell+1}=(1-\gb)^2 \gD_\ell +0.01\ga_{_{\ell+1}}
 (\th_{\ell+1}^3 n)^{1/2} $, \raf{v3bb} and $\ga_j \ll (\th_j^3 n)^{1/2}$ give
  \bean P_{\ell+1} &\leq & \pr \bb[ |B_\ell - b_{\ell}| >\sfrac{1}{400}
  \ga_{_{\ell+1}}   (\th_{\ell+1}^3
n)^{1/2}\bb]  \\ & & + \pr\bb[ |M_\ell - (1-(1-\gb)^2) N_\ell|>
\sfrac{1}{200}\ga_{_{\ell+1}}   (\th_{\ell+1}^3 n)^{1/2},
|N_{\ell} -
n_{_{\ell}} | \leq   \gD_{\ell}\bb]\\
&\leq & e^{-\Omega(\ga_{\ell+1}^2)}+ \pr\bb[ |M_\ell -
(1-(1-\gb)^2) N_\ell|> \sfrac{1}{200}\ga_{_{\ell+1}}
(\th_{\ell+1}^3 n)^{1/2}\bb| |N_{\ell} - n_{_{\ell}} |\leq
\gD_{\ell} \bb]. \eean The desired bound follows, since \raf{v3mm}
yields
$$
\pr\bb[ |M_\ell - (1-(1-\gb)^2) N_\ell|>
\sfrac{1}{200}\ga_{_{\ell+1}} (\th_{\ell+1}^3 n)^{1/2}\bb|
N_{\ell}  \bb] \leq e^{-\Omega(\ga_{\ell+1}^2)}$$
 for given
$N_\ell$ with $|N_{\ell} - n_{_{\ell}} |\leq \gD_{\ell} \ll
\th_{\ell}^3 n$.

 The second inequality  holds for  \raf{v4hh}
 gives
 $$
\pr \bb[ \exists ~j=1,...,\ell ~~s.t.~~ |H_j -2\th_j \gl
(1-e^{-\th_j \gl}-\th_j\gl
 e^{-2\th_j\gl})n | \geq \gD_j  \, \bb]
 \leq 2\sum_{j=1}^\ell e^{-\Omega(\ga_j^2)}
 =e^{-\Omega(\ga_{\ell}^2)}. $$
 \qed

\old{

 Suppose $\gl = 1- \gs $ with $\gs \gg n^{-1/3}$.
  Then we take $\ga= \log (\gs^3
n)$ and $\gb$ in the range $1/3\leq \gb \leq 1/2$ so that $
(1-\gb)^a = n^{-1/3} \log (\gs^3 n)$ for an integer $a$. As the
ratio $(1-\gb)^{\frac{j}{4} } \th_j^{-3/2} = (1-\gb)^{-5j/4}$ has
the maximum when $j=a$ and $\sum_{i=1}^{j} (1-\gb)^{\frac{j-i}{4}}
= O(1/\gb)$, \beq{v4aap} \ga \sum_{i=1}^{j}
(1-\gb)^{\frac{2j-i-a}{4}}= \ga(1-\gb)^{\frac{j-a}{4} }
\sum_{i=1}^{j} (1-\gb)^{\frac{j-i}{4}} \ll  (\th_j^3 n)^{1/2},\eeq
for $j=1, ..., a$.

\begin{lem} \label{v3sub}  Suppose $\gl = 1- \gs $ with $\gs \gg
n^{-1/3}$,  $\ga= \log(\gs^3 n)$ and  $\gD_j=0.1 \ga (
(\gs+\th_j)\th_j^2 n)^{1/2} \sum_{i=1}^{j}
(1-\gb)^{\frac{2j-i-a}{4}}. $ Then, for $\ell=1, ...,a$,
$$ \pr \bb[ \, \exists ~j=1,...,\ell ~~s.t.~~
 |N_{j} - 2 \th_j^2\gl  ( 1-\gl e^{-2\th_j \gl })n| >  \gD_{j}   \bb] \leq
  e^{-\Omega(\ga^2(1-\gb)^{\frac{\ell-a}{2}})},
$$
and
$$ \pr \bb[ \, \exists ~j=1,...,\ell ~~s.t.~~
 |L_{j} - 2 \th_j\gl  ( \th_j- 1+ e^{-\th_j \gl })n| >  2\gD_{j}   \bb] \leq
   e^{-\Omega(\ga^2(1-\gb)^{\frac{\ell-a}{2}})}.$$
\end{lem}

\noindent \pf  Let $n_j = 2 \th_j^2\gl  ( 1-\gl e^{-2\th_j \gl
})n$ and $\ga_j = (1-\gb)^{\frac{j-a}{4}} \ga$. Then \raf{v4aap}
gives $\ga_j \ll (\th_j^3 n)^{1/2}$.
  For the first inequality, it is enough to show that
$$ \pr \bb[ |N_{j+1} - n_{j+1}| >  \gD_{j+1} \bb|
|N_{j} - n_{j}| \leq \gD_j \bb] \leq e^{-\Omega(\ga_{j+1}^2)},
$$
as $\sum_{j=1}^\ell e^{-\Omega(\ga_{j}^2)} =
e^{-\Omega(\ga_\ell^2)}$.

Since $N_{j+1} = N_j - M_j -2B_j$, $\ga_{j+1} \ll (\th_{j+1}^3
n)^{1/2}, $ and $n_j =\Theta((\gs+\th_j)\th_j^2 n)$,  \raf{v3mm}
and \raf{v3bb} give
$$ \pr \bb[ |N_{j}- M_j  -  (1-\gb)^2 N_j  |
\geq \sfrac{1}{20}\ga_{j+1}   (\th_{j+1}^3  n)^{1/2}\bb| |N_{j} -
n_{j}| \leq \gD_j \bb]\leq e^{-\Omega(\ga_{j+1}^2)},
$$ and
$$ \pr \bb[ |B_j -(\th_{j+1}\gl)^2 e^{-2\th_{j+1} \gl}
(1- e^{-2\gb \th_j \gl})n| \geq \sfrac{1}{40}\ga_{j+1}   (\th_j^3
n)^{1/2} \bb] \leq e^{-\Omega(\ga_{j+1}^2)}. $$ Hence, using
$N_{j+1} = N_j -M_j - 2B_j$ and  $n_{j+1} = (1-\gb)^2 n_j -
2(\th_{j+1}\gl)^2 e^{-2\th_{j+1} \gl} (1- e^{-2\gb \th_j \gl})n$,
$$ \pr \bb[ |N_{j+1} - n_{j+1} | \geq (1-\gb)^2 \gD_j + 0.1\ga_{j+1}
 (\th_{j+1}^3 n)^{1/2} \bb| |N_{j} - n_{j}| \leq
\gD_j \bb] \leq e^{-\Omega(\ga_{j+1}^2)}. $$ The desired
inequality follows by  $\gD_{j+1}\geq (1-\gb)^2 \gD_j +0.1
\ga_{j+1}
 ((\gs+\th_j)\th_{j+1}^2 n)^{1/2} $.

 The second inequality also  follows since \raf{v4hh}
 gives
 $$
\pr \bb[ \exists ~j=1,...,\ell ~~s.t.~~ |H_j -2\th_j \gl
(1-e^{-\th_j \gl}-\th_j\gl
 e^{-2\th_j\gl})n | \geq \gD_j  \, \bb]
 \leq 2\sum_{j=1}^\ell e^{-\Omega(\ga_j^2)}
 =e^{-\Omega(\ga_\ell^2)}. $$
 \qed

} 

\mn

\noindent {\bf Proof of Main Lemma.}  We first estimate the
probability that all light clones
 disappear during phase $j$. Observe that the number of
 light clones  is bounded by renewal random
walk processes  with  negative drift: If the chosen light clone is
matched to another light clone, then the number decreases by $2$.
If it is matched to a clone of a variable with type larger than or
equal to $(2,2)$, the number decreases by $1$. If it is matched to
a clone of  a variable with type $(1,b)$ or $(b,1)$, $b\geq 2$,
then the number decreases, in expectation, by $-1+ \frac{b}{b+1}$.
Thus, there is absolute constant $h>0$ such that the expected
number of light clones is less than $-h$.

If all light clones
 disappear during phase $j$, $j=1,..., a-1$, then, either $L_{j+1}\leq 0.1
 \ga (\th_{j+1}^3 n)^{1/2} $ or the  renewal random walks
with negative drift must reach beyond $0.1
 \ga (\th_{j+1}^3 n)^{1/2} $. Lemma \ref{v3sup} gives the
 probability of the former is
 $e^{-\Omega(\ga^2(1-\gb)^{\frac{\ell-a}{2}})}$.
 For the latter, observe that
 the total number of walks is  less than $N_j$, which is $O(\th_j^3 n)$
with probability  $e^{-\Omega(\ga^2 (1-\gb)^{\frac{j-a}{2}})}$. We
consider  excursions that are  segments of
 the renewal random walks between two consecutive visits to  $0$.
The generalized Chernoff bound, or  Corollary \ref{nd}, with
$a_{i}, b_i, \vp, \gd =\Theta (1)$ yields that
 the probability of each  excursion reaching beyond
 $0.1 \ga (\th_{j+1}^3 n)^{1/2} $ is at most
$$
\sum_{m\geq 1} e^{-\Omega(\min \{\ga (\th_{j+1}^3 n)^{1/2}+hm   ,
\, \, (\ga (\th_{j+1}^3 n)^{1/2}+hm)^2 / m\} } \leq \sum_{m\geq 1}
e^{-\Omega(\ga (\th_{j+1}^3 n)^{1/2})-\Omega(m) } = e^{-\Omega(\ga
(\th_{j+1}^3 n)^{1/2})} .
$$
As there are at most $N_j$ excursions
such an excursion exists with probability at most
$$ e^{-\Omega(\ga^2(1-\gb)^{\frac{j-a}{2}})}
+ O\bb(\th_j^3 n  e^{-\Omega(\ga (\th_{j+1}^3
 n)^{1/2})}\bb)\leq
 e^{-\Omega(\ga^2(1-\gb)^{\frac{j-a}{2}})} +
 e^{-\Omega(\ga^2(1-\gb)^{\frac{3(j-a)}{2}})}
 = e^{-\Omega(\ga^2(1-\gb)^{\frac{j-a}{2}})} . $$

Therefore, there exists no light clones in a step of the $j^{\rm
th}$ phase is at most
$e^{-\Omega(\ga^2(1-\gb)^{\frac{\ell-a}{2}})} $, which yields
$$ \pr [ \gL_C \geq \thl\gl + \ga (\thl n)^{-1/2}] \leq
e^{-\Omega(\ga^2)}, $$ replacing $\ga$ by $\ga/\gl$.

 On the other hand, after
$(a-1)^{\rm th} $ phase, there are at most  $O(\ga (\thl^3
n)^{1/2})$ light clones and $\Omega(\thl^3 n)$ unmatched clones of
variables of type larger than $(1,1)$, with probability
$1-e^{-\Omega(\ga^2)}$. As the number of light clones has negative
drift, no light clone exists after $ O(\ga (\thl^3
n)^{1/2})=O(\frac{\ga}{(\thl^3 n)^{1/2}} \thl^3 n)$ more clones
are matched with probability $1-e^{-\Omega(\ga^2)}$. Thus, the
cut-off value when no light clone exists for the first time cannot
be smaller than $(1-O(\frac{\ga}{(\thl^3 n)^{1/2}}) )\thl \gl =
\thl\gl - O(\ga (\thl n)^{-1/2})$  with probability
$1-e^{-\Omega(\ga^2)}$ by the cut-off line lemma, or \raf{v3mm}.
Replacing $\ga$ by $c\ga$ for appropriate constant, we conclude
that
$$ \pr [ \gL_C \leq \thl\gl - \ga (\thl n)^{-1/2} ]\leq e^{-\Omega(\ga^2)}. $$
\qed

\old{
 The same analysis together with \raf{v4hh} and Lemma
\ref{v3sub} also yields the following lemma.
\begin{thm}
\label{lower} If $\gl =1 -\gs$ is uniformly bounded from below  by
$0$ with $\gs \gg n^{-1/3}$, then
$$ \pr [ \gL_C \geq  n^{-1/3} \log (\gs^3 n)] \leq e^{-\Omega(\log^2
(\gs^3n))},$$ and for the number $L$ of light clones when the
cut-off line reaches $n^{-1/3} \log (\gs^3 n)$
$$ \pr \bb[ L - 2\gs n^{1/3} \gl\log^2 (\gs^3 n) \leq (\gs n^{1/3})^{1/2}
\log^2(\gs^3 n)
 \bb] =e^{-\Omega(\log^2 (\gs^3n))}.$$
\end{thm}
}

\bn

  Let $F^*_{C}(n,p)$ be the core of
$F_{PC}(n,p)$. One may define $C^*_{n,p}$, $C^*_{n,p}(i,j)$,
$D^*_{n,p} ( i,j)$ and $M^*_{n,p}( i,j)$ for $F_{PC} (n,p)$  as
$C_{n,p}$, $C_{n,p}(i,j)$, $D_{n,p} ( i,j)$ and $M_{n,p}( i,j)$
are defined for $F (n,p)$. We first
 estimate  $|D^*_{n,p} (
i,j)|$ and
 $M^*_{n,p}(
i,j)$. Notice that the upper and lower
 bounds for $|D^*_{n,p} (
i',j')|$'s yield bounds for
 $|C^*_{n,p}(i,j)|$ as
 $$ 
   |C^*_{n,p}(i,j)|=|D^*_{n,p}( i,j)|-  |D^*_{n,p}( i+1,j)|
 -|D^*_{n,p}( i,j+1)|+|D^*_{n,p}(i+1,j+1)|.
$$ 
Let $D^{\pm}_\ga (i,j)$ be the sets of variables that are of
type larger than or equal to $(i,j)$ at $\mu^{\pm}_\ga:=\thl\gl\pm
\ga(\thl n)^{-1/2}$, respectively, and let $ M_\ga^{\pm} (i,j)$ be
the number of clones of variables in $D_\ga^{\pm} (i,j)$ less than
$\mu_\ga^{\pm}$, respectively. Then,  Lemma \ref{upper} gives
$$ \pr \bb[D_\ga^{-}(i,j) \sub D^*_{n,p}( i,j)  \sub D_\ga^{+}(i,j)
~~\mbox{for all $i,j$ }\bb] =1-e^{-\Omega(\ga^2)}, $$ and
$$ \pr \bb[M_\ga^{-}(i,j) \leq   M^*_{n,p} (i,j)  \leq  M_\ga^{+}(i,j)
~~\mbox{for all $i,j$ }\bb] =1-e^{-\Omega(\ga^2)}. $$ Since
$|D_\ga^{\pm}(i,j)|$ and $M_\ga^{\pm}(i,j)$ are the sums of i.i.d
random variables and it is easy to check all the conditions of the
generalized Chernoff bound with $a_i, b_i = \Theta (\thl^{i+j}),
~\xi_{_0}=\gd =1 $, we have
  $$\pr \bb[\, \, \bb| |D_\ga^{\pm}(i,j)|  -
  Q_i(\mu^{\pm}_\ga) Q_j(\mu^{\pm}_\ga)n \bb|\geq
\gD\, \bb] \leq 2e^{-\Omega(\min\{\gD, \frac{\gD^2}{\thl^{i+j
}n}\})}, $$ respectively, and
$$\pr \bb[\, \, \bb|M_\ga^{\pm}(i,j)  - \mu^{\pm}_\ga \bb(  Q_{i-1}( \mu^{\pm}_\ga  )Q_j(
\mu^{\pm}_\ga ) +\ Q_i( \mu^{\pm}_\ga )Q_{j-1}( \mu^{\pm}_\ga
)\bb)n
 \bb|\geq \gD\, \bb] \leq
 2e^{-\Omega(\min\{\gD, \frac{\gD^2}{\thl^{i+j}n}\})}, $$
 respectively.
Therefore,
$$\pr\bb[ |D^*_{n,p} (i,j)| -  Q_i (\mu_\ga^+) Q_j (\mu_\ga^+) n\geq  \gD
\bb] \leq 2e^{-\Omega(\min\{\gD, \frac{\gD^2}{\thl^{i+j}n}\})} +
e^{-\Omega(\ga^2)} ,$$
$$ \pr\bb[ |D^*_{n,p} (i,j)|  - Q_i (\mu_\ga^-) Q_j (\mu_\ga^-) n\leq
-\gD \bb] \leq 2e^{-\Omega(\min\{\gD,
\frac{\gD^2}{\thl^{i+j}n}\})}+e^{-\Omega(\ga^2)},
$$
and
$$ \pr \bb[ M_{n,p}^* (i,j)
-\bb(   Q_{i-1}( \mu^{+}_\ga  )Q_j( \mu^{+}_\ga ) + Q_i(
\mu^{+}_\ga )Q_{j-1}( \mu^{+}_\ga )\bb) \mu^{+}_\ga n \geq \gD\bb]
\leq 2e^{-\Omega(\min\{\gD,
\frac{\gD^2}{\thl^{i+j}n}\})}+e^{-\Omega(\ga^2)}, $$
$$
\pr \bb[ M_{n,p}^* (i,j) - \bb( Q_{i-1}( \mu^{-}_\ga )Q_j(
\mu^{-}_\ga ) + Q_i( \mu^{-}_\ga )Q_{j-1}( \mu^{-}_\ga )\bb)
\mu^{-}_\ga n \leq -\gD\bb] \leq 2e^{-\Omega(\min\{\gD,
\frac{\gD^2}{\thl^{i+j}n}\})}+e^{-\Omega(\ga^2)}.$$
  We also have
 \begin{eqnarray}\label{ebr2}   \pr [ |D^*_{n,p} (i,j)| \geq \ell
 ]&\leq  &e^{-\Omega{(\ga^2)}} + \pr[ |D^+_{n,p} (i,j)| \geq
 \ell]  \nonumber \\
 & \leq &
e^{-\Omega{(\ga^2)}} +  {n \choose \ell}( Q_i (\mu_\ga^+) Q_j
(\mu_\ga^+) )^{\ell} \nonumber  \\
 & \leq &
e^{-\Omega{(\ga^2)}} + \frac{(Q_i (\mu_\ga^+) Q_j (\mu_\ga^+)
n)^{\ell}}{\ell !}.  \end{eqnarray}

If $(i,j)$ is  fixed, it is easy to see that
  $$Q_i (\mu_\ga^\pm) Q_j (\mu_\ga^\pm)= Q_i(\thl\gl)Q_j (\thl\gl)
  +O(\ga \thl^{i+j-1}  (\thl n)^{-1/2}),$$
and similarly
  $$ \mu^{\pm}_\ga Q_i( \mu^{\pm}_\ga )Q_{j}(
\mu^{\pm}_\ga ) =   \thl\gl Q_i( \thl\gl)Q_{j}( \thl\gl )+O(\ga
\thl^{i+j} (\thl n)^{-1/2}).
$$
Replacing $\ga$ by $c\ga$ for an appropriate constant $c>0$ and
taking $\gD= c\ga \thl^{i+j-1}(\thl n)^{-1/2} n= c \ga
\thl^{i+j-3} (\thl^3 n)^{1/2} $, we have, for $Q_i = Q_i
(\thl\gl)$,
$$\pr\bb[\,\, \bb| |D^*_{n,p} (i,j)| -  Q_i  Q_j  n
\bb| \geq   \ga \thl^{i+j-3} (\thl^3 n)^{1/2}  \bb] \leq
e^{-\Omega(\thl^{i+j-3} \ga^2)}+ e^{-\Omega(\ga^2)},$$ and
$$
\pr\bb[\,\, \bb|M^*_{n,p} (i,j) -  \thl\gl \bb(Q_{i-1} Q_j + Q_{i}
Q_{j-1} \bb) n \bb| \geq \ga \thl^{i+j-3} (\thl^3 n)^{1/2} \bb]
\leq e^{-\Omega(\thl^{i+j-3} \ga^2)}+ e^{-\Omega(\ga^2)}.$$

Thus, for $F(n,p)$, Theorem \ref{equiv} gives

\begin{thm}
\label{main} Suppose  $p(2n-1) =1+\gs $ is uniformly bounded from
above with $\gs \gg n^{-1/3}$. Let $\gl=1+\gs$, $\gD>0$, $1\ll
\ga\ll (\thl n)^{1/2} $, and $\mu^{\pm}_\ga = \thl \gl \pm \ga
(\thl n)^{-1/2}$, respectively. Then,
$$\pr\bb[ |D_{n,p} (i,j)| -  Q_i (\mu_\ga^+) Q_j (\mu_\ga^+) n\geq  \gD
\bb] \leq 2e^{-\Omega(\min\{\gD,
\frac{\gD^2}{\thl^{i+j}n}\})}+e^{-\Omega(\ga^2)},$$
$$ \pr\bb[ |D_{n,p} (i,j)|  - Q_i (\mu_\ga^-) Q_j (\mu_\ga^-) n\leq
\gD \bb] \leq 2e^{-\Omega(\min\{\gD,
\frac{\gD^2}{\thl^{i+j}n}\})}+e^{-\Omega(\ga^2)},
$$
and
$$ \pr [ |D_{n,p} (i,j) | \geq \ell ] \leq O
\bb(\frac{((1+\frac{\ga}{(\thl n)^{1/2}})\thl
\gl)^{(i+j)\ell/2}}{(\ell!)^{1/2}}\bb) + e^{-\Omega (\ga^2)}.
$$
We also have
$$ \pr \bb[ M_{n,p} (i,j)
-\mu^{+}_\ga\bb(   Q_{i-1}( \mu^{+}_\ga  )Q_j( \mu^{+}_\ga ) +
Q_i( \mu^{+}_\ga )Q_{j-1}( \mu^{+}_\ga )\bb) n \geq \gD\bb] \leq
2e^{-\Omega(\min\{\gD,
\frac{\gD^2}{\thl^{i+j}n}\})}+e^{-\Omega(\ga^2)}, $$
$$
\pr \bb[ M_{n,p} (i,j) -\mu^{-}_\ga \bb( Q_{i-1}( \mu^{-}_\ga
)Q_j( \mu^{-}_\ga ) + Q_i( \mu^{-}_\ga )Q_{j-1}( \mu^{-}_\ga )\bb)
n \leq \gD\bb] \leq 2e^{-\Omega(\min\{\gD,
\frac{\gD^2}{\thl^{i+j}n}\})}+e^{-\Omega(\ga^2)}.$$ In particular,
for fixed $(i,j)$,
$$\pr\bb[\,\, \bb| |D_{n,p} (i,j)| -  Q_i  Q_j  n
\bb| \geq   \ga \thl^{i+j-3} (\thl^3 n)^{1/2}  \bb] \leq
e^{-\Omega(\thl^{i+j-3} \ga^2)}+ e^{-\Omega(\ga^2)},$$ and
$$
\pr\bb[\,\, \bb|M_{n,p} (i,j) -  \thl\gl \bb(Q_{i-1}  Q_j + Q_{i}
Q_{j-1} \bb) n \bb| \geq \ga \thl^{i+j-3} (\thl^3 n)^{1/2} \bb]
\leq e^{-\Omega(\thl^{i+j-3} \ga^2)}+ e^{-\Omega(\ga^2)}.$$
\end{thm}
\qed

Theorem \ref{main} together with
$$
   |C _{n,p}(i,j)|=|D _{n,p}( i,j)|-  |D _{n,p}( i+1,j)|
 -|D _{n,p}( i,j+1)|+|D _{n,p}(i+1,j+1)|.
$$
 implies that
$$\pr\bb[\,\, \bb| |C _{n,p} (i,j)| -  P_i P_j  n
\bb| \geq   \ga \thl^{i+j-3} (\thl^3 n)^{1/2}  \bb] \leq
e^{-\Omega(\thl^{i+j-3} \ga^2)}+ e^{-\Omega(\ga^2)}.$$ Similarly,
if $i\geq j$, then
$$ |D _{n,p} (i,j)\cup D _{n,p} (j,i)|=
|D _{n,p} (i,j)| +| D _{n,p} (j,i)|-|D _{n,p} (i,i)|$$ gives
$$\pr\bb[\,\, \bb| |D _{n,p} (i,j)\cup D _{n,p} (j,i)|
- \bb( 2 Q_i Q_j  n -  Q_i Q_i  \bb) n\bb| \geq \ga \thl^{i+j-3}
(\thl^3 n)^{1/2} \bb] \leq e^{-\Omega(\thl^{i+j-3} \ga^2)}+
e^{-\Omega(\ga^2)}.$$ The  last two bounds of Theorem \ref{seccor}
are already in Theorem \ref{main}.

Furthermore,  as $C_{n,p} = D_{n,p} (1,1)$, $K_{n,p}= D_{n,p}
(2,1) \cup D_{n,p} (1,2)$, $2|F_C (n,p)|= M_{n,p} (1,1)$, $2|F_K
(n,p)| = M_{n,p} (1,2)+M_{n,p}(2,1)- M_{n,p} (2,2)$ and $Q_1
(\thl) = 1-e^{-\thl\gl} =\thl$, Corollary \ref{corekernel} follows
from Theorem \ref{seccor}.

Finally, the first two bounds in Corollary \ref{swthm} follow from
Theorem \ref{seccor} since $\thl= \Theta(\gs)$. For the last
bound, if $\gs^4 n \leq \eps$ for a small positive constant
$\eps$, then Theorem \ref{seccor} implies that all variables in
the core are of types $(1,1)$, $(1,2)$ or $(2,1)$, with
probability $1- O(\eps^2)$. If $\gs^4 n \geq \gb $ for a large
constant $\gb>0$, then Theorem \ref{seccor} with $\ga = \gb^{-0.1}
(\thl^3 n)^{1/2}$ also gives
$$ \pr [ M_{n,p} (2,2) - 2\thl\gl Q_1(\thl\gl)Q_2(\thl\gl) n \geq \gb^{-0.1}
\thl^4 n ] \leq e^{-\Omega(\gb^{-0.2} \thl^4 n )} \leq
e^{-\Omega(\gb^{0.8} )} . $$ Similar bounds hold for $M_{n,p}
(1,3)$ and $M_{n,p} (3,1)$. As  $\thl= \Theta (\gs)$ and $\thl\gl
Q_1(\thl\gl)Q_2(\thl\gl)= \Theta (\thl^4)$, the desired bound
follows.  If $\eps \leq \gs^4 n \leq \gb$, we simply use the bound
for $\gl'= 1+(\gb/n)^{1/4}$, or $p'=
\frac{1+(\gb/n)^{1/4}}{2n-1}$: Since $\pr [ M_{n,p} \geq h] \leq
\pr [ M_{n,p'} \geq h]$ and $M_{n,p'}= O(\gb )=O((\gb/\eps) \gs^4
n) $ with probability $1-e^{-\Omega(\gb^{0.8})}$, the bound
follows.

\old{
 Moreover, for $(i,j)\geq (2,2)$, any $\gd$ in the range
$\thl^{i+j-3} <\gd \ll  1$  and $\ga= (\thl^3 n)^{1/2}$,
$$\pr \bb[ |D _{n,p} (i,j)| \geq \sfrac{(\thl\gl)^{i+j}}{i!j!}
+ \gd \ga (\thl^3 n)^{1/2}\bb] \leq e^{-\Omega(\min\{\gd \ga
(\thl^3 n)^{1/2}, \frac{ \gd^2 \ga^2}{\thl^{i+j-3}}\})}+
e^{-\Omega(\gd \ga^2)}= e^{-\Omega(\gd \ga^2)},
$$
and
$$\pr \bb[ M _{n,p} (i,j) \geq \sfrac{(\thl\gl)^{i+j}}{(i-1)!j!}
+ \sfrac{(\thl\gl)^{i+j}}{i!(j-1)!} + \gd \ga (\thl^3 n)^{1/2}\bb]
\leq e^{-\Omega(\min\{\gd \ga (\thl^3 n)^{1/2}, \frac{ \gd^2
\ga^2}{\thl^{i+j-3}}\})}+e^{-\Omega(\gd \ga^2)}= e^{-\Omega(\gd
\ga^2)}.
$$} 

\section{Scaling Window:  Proofs of Theorems \ref{config}, \ref{subker} and
\ref{sub2sat}} \label{sw2}

 Suppose  all $C_{n,p}(i,j)$'s  for
$(i,j)\geq (1,1)$ are given. The first thing we need to establish
is that all $2$-SAT formulae with the same $C_{n,p}(i,j)$'s are
equally likely to be the core of $F(n,p)$. More generally, it is
not hard to show the following lemma.

\begin{thm} (Restated)
Suppose  two formulae 
have the same number of clauses on the same number of underlying
variables, and all underlying variables are of type at least
$(1,1)$. Then the two are equally likely to be the core of
$F(n,p)$.
\end{thm}

\noindent \pf Let $F_1$ and $F_2$ be the two formulae. After an
appropriate permutation, we may assume that  the two formulae have
the same set of underlying variables. Then, a formula having $F_1$
as its core can be mapped to the formula obtained by replacing
 clauses in $F_1$ with clauses of $F_2$. It is easy to see
that the core of the formula obtained this way is $F_2$. It is
also clear that the map is one-to-one and onto. Furthermore, two
formulae mapped each other have the same number of clauses, which
means that the random formula $F(n,p)$ is equally likely to be one
of  the two formulae. \qed

\mn

\newcommand{\prc}{\pr\,\! \! ^*}
 We now  consider the configuration model for given $C_{n,p}
(i,j)$: Similar to the Poisson cloning model, take $i$ clones of
$x$ and $j$ clones of $\bar{x} $ for each variable $x\in C_{n,p}
(i,j) $. The uniform random perfect matching on all clones is
called the random configuration. The random configuration then
yields a $2$-SAT multiformula after contractions. The event that
the multiformula has neither loops nor multiple clauses is called
$SIMPLE$ or $SIM$. Conditioned on $SIM$, the random $2$-SAT
formula has the uniform distribution among all $2$-SAT formulae
with the same $C_{n,p}(i,j)$'s. This is not difficult to see as
the number of perfect matchings that yield a fixed $2$-SAT formula
 is $\prod_{(i,j)} (i!j!)^{|C_{n,p}(i,j)|}$.
It is known 
that the probability of $SIM$ is uniformly  bounded below from
$0$, especially, for any event $A$ in the uniform model, or
equivalently in the configuration model,
  \beq{sim} \pr [A]= \prc [ A|SIM] \leq \prc [SIM]^{-1} \prc[A]
  = O(\prc [A]), \eeq
  where the probability $\prc$ is taken over the random configuration
   without
  any condition.
Hence,  as far as the constant factor is not concerned, it is
enough to bound the desired probabilities in the configuration
without any condition. To clarify terminology, we recall that the
configuration model is obtained from the random configuration  by
conditioning $SIM$. For an event $A$ depending on the random
configuration only, such as the event that the $i^{\rm th}$ clone
of $y$ and the  $j^{\rm th}$ clone of $z$ are matched, $\pr[A]$
may not be well-defined, but $\prc[A]$ or $\pr^* [A|SIM]$ may be
still considered.

We may be able to estimate the probability of $SIM$ in the case
that all but few clones  are clones of type $(1,1)$ literals:
Suppose all $N$ but $o(N^{1/2})$ clones are clones of type $(1,1)$
literals. First, with probability  $1-o(1)$,  no pair of clones
that are not clones of type $(1,1)$ literals is matched. Thus, the
multiformula is not simple mainly because two clones of type
$(1,1)$ variable $x$ and its negation are matched. Let $A_x $ be
such an event. Then, for the set $U$ of all type $(1,1)$
variables,
$$ \sum_{\ell=0}^{2i+1} (-1)^\ell \sum_{W\sub U \atop |W|=\ell}
\prc \bb[ \bigcap_{x\in W}  A_x \bb] \leq  \prc \bb[ \, \,
\overline{\bigcup_{x\in W} A_x}  \, \,\bb] \leq \sum_{\ell=0}^{2i}
(-1)^\ell \sum_{W\sub U \atop |W|=\ell} \prc \bb[ \bigcap_{x\in W}
  A_x \big] ,$$ for all $i\geq 0$.
For $\ell=o(N^{-1/2} ) $ and $|W|=\ell$,
$$ {|U| \choose \ell} ={ N/2-o(N^{1/2} ) \choose \ell } =
\frac{(1+O(\frac{\ell(\ell+
N^{1/2})}{N}))(\frac{N}{2})^\ell}{\ell!} $$ and $$\prc [
\cap_{x\in W}
  A_x] = \frac{(N-2\ell-1)!!}{(N-1)!!}
  = \frac{1+O(\frac{\ell^2}{N})}{N^{\ell}}. $$
Therefore,
 \beq{esim} \prc \bb[\,\, \overline{\bigcup_{x\in W} A_x}\, \, \bb] = (1+o(1)) e^{-1/2},
~~{\rm and} ~~ \prc [ SIM] = (1+o(1)) e^{-1/2} \eeq

We are now ready to prove Theorems \ref{subker} and \ref{sub2sat}.

 \noindent {\bf Proof of Theorem \ref{subker}} We
may generate $F(n,p)$ with $\gl_p:= p(2n-1)= 1-\gs$ by first
taking $F(n, q)$ with $\gl_q= 1+n^{-1/3}\log (\gs^3 n)$ and then
independently selecting each clause of $F(n, q)$ with probability
$p/q= \frac{1-\gs}{1+n^{-1/3} \log (\gs^3 n)}=1-\gs +o(\gs)$.

Applying  Corollary \ref{corekernel} for $\gl_q$ and $\ga =
\log(\gs^3 n)$ and using  $\th_q:=\th_{\gl_q} = 2n^{-1/3} \log
(\gs^3 n)+ O(n^{-2/3} \log^2 (\gs^3 n)) $, we have
$$ |C_{n,q} | =  \th_{q}^2 q^2 n +O(\th_{q}^3  n)+ O(
(\th_{q} n)^{1/2}\log (\gs^3n))= (4+o(1))n^{1/3} \log^2 (\gs^3 n),
$$
 $$ |K_{n,q} |= \th_{q}^3 \gl^3 n + O(\th_{q}^4  n)+
O( (\th_{q}^3 n)^{1/2}\log (\gs^3n))= (8+ o(1))\log^3 (\gs^3 n)
$$
and
$$ |F_C (n,q)| = \th_{q}^2 \gl^2 n +O(\th_{q}^3 n)
 +O((\th_{q} n)^{1/2}\log (\gs^3n)) = (4+ o(1))n^{1/3} \log^2
 ( \gs^3 n),
  $$
  $$ |F_K (n,q)| = \sfrac{3}{2}\th_{q}^3 \gl^3 n +O(\th_{q}^4 n)
 +O((\th^3_{q} n)^{1/2}\log (\gs^3n)) = (12+ o(1))
 \log^3
 ( \gs^3 n),
  $$ and for $K_{n,q} (1,2) := C_{n,q} (1,2) \cup C_{n,q} (2,1)$
$$ |K_{n,q} (1,2)| = (8+ o(1))\log^3 (\gs^3 n)$$
  with probability $1-e^{-\Omega (\log^2 (\gs^3 n))}$.
Theorem \ref{main} with the same $\ga$ also gives
 \beq{nolarge} |D_{n,q} (i,j)]=0~~\mbox{if ~$i+j \geq 7$} \eeq with
probability $1-n^{-7/6+o(1)}$.

 Suppose
$C_{n,q} (i,j)$'s are given with $C_{n,q} :=\cup_{(i,j)\geq (1,1)}
C_{n,q} (i,j)$, $K_{n,q} :=\cup_{(i,j)> (1,1)} C_{n,q} (i,j)$,
$D_{n,q}(i,j) :=\cup_{(i',j')> (i,j)} C_{n,q} (i',j')$, and
$K_{n,q} (1,2) := C_{n,q} (1,2) \cup C_{n,q} (2,1)$ satisfying the
above conditions, and the number $M_C$ and $M_K$ of clones of
variables in $C_{n,q}$ and $K_{n,p}$, respectively, satisfy
$$M_C =\sum_{(i,j)\geq (1,1) }(i+j)| C_{n,q}
(i,j)| = (8+o(1)) n^{1/3} \log^2 (\gs^3 n),$$ and   $$ M_K
=\sum_{(i,j)> (1,1) }(i+j)| C_{n,q} (i,j)| = (24+o(1)) \log^3
(\gs^3 n). $$
 In the
corresponding random configuration  for given $C_{n,q} (i,j)$'s,
two clones $y,z$ of variables in $K_{n,p}$ may yield a clause
after resolutions of variables in $C_{n,q}(1,1)$. This occurs if
and only if there are $x_{_1}, ..., x_{_\ell} \in C_{n,q}(1,1)$
such that $ \{w_{_0}:=y, w_1\}, \{\bar{w}_1, w_2\} , ...,
\{\bar{w}_{\ell-1}, w_\ell\}, \{\bar{w}_{\ell}, w_{\ell+1}:=z\} $
are edges in the random configuration, where $w_i, \bar{w}_i$ are
the two clones of $x_i$ and $\bar{x}_i$ (not necessarily
respectively), including the case $\ell=0$. If this event occurs,
we say that the length $\ell(y,z)$ of $y,z$ is $\ell+1$ and the
edges $\{w_{i}, w_{i+1}\}$ (resp. the corresponding clauses after
contractions) are called intermediate edges (resp. clauses) of the
pair.   The length $\ell(y,z) $ is infinity if no such $x_i$'s
exist.  It is easy to see that $\prc [ \ell (y,z) = 1]
=\frac{1}{M_C-1}.$ Similarly,
$$\prc [ \ell (y,z) = 2]
=\bb(1-\frac{M_K-1}{M_C-1}\bb)\frac{1}{M_C-3}, $$ and, in general,
$$\prc [ \ell (y,z) = \ell]
=\bb(1-\frac{M_K-1}{M_C-1}\bb)\bb(1-\frac{M_K-1}{M_C-3}\bb) \cdots
\bb(1-\frac{M_K-1}{M_C-2\ell+3}\bb) \frac{1}{M_C-2\ell+1} .$$ For
$i=1,...,4$ and $\ell_1, ..., \ell_i \leq \frac{10}{\gs}\log(\gs^
3n) \ll n^{1/3}$, the same argument also gives
 \begin{eqnarray}
 \prc [ \ell (y_{_j},z_{_j} ) = \ell_{j}, ~~j=1,...,i]
&=& \bb(1- \frac{(1+o(1))M_K}{M_C}\bb)^{\ell_{1}+\cdots+\ell_{i}}
\bb(\frac{1+o(1)}{M_C }\bb)^i  \nonumber \\ &=&
 \bb(\frac{ 1+o(1)}{8n^{1/3} \log^2 (\gs^3 n)}\bb)^i .
 \end{eqnarray}

Notice that a pair $y,z$ of clones of variables in $K_{n,q}$
yields the corresponding clause in the kernel $F_K(n,p) $
 only if $\ell(y,z) < \infty$ and all the  $\ell(y,z)$
intermediate clauses of the pair are in $F(n,p)$.
 Such an event occurs with probability
 \bean \sum_{\ell \geq 1}
(p/q)^{\ell}\prc[ \ell(y,z) =\ell|SIM] & =& O\bb(  \sum_{\ell =
1}^{\frac{10}{\gs}\log(\gs^ 3n)} (1-\gs)^{\ell} \prc[ \ell(y,z)
=\ell] 
+ (1-\gs)^{\frac{10}{\gs}\log(\gs^ 3n)}\bb)
\\
&=& O\bb(\sum_{\ell = 1}^{\frac{10}{\gs}\log(\gs^ 3n)}
(1-\gs)^{\ell} \prc[ \ell(y,z) =\ell] \bb) +O( (\gs^3 n)^{-10})
 . \eean

Similarly, if   the kernel $F_{K} (n,p)$ of $F(n,p)$ has $i$ or
more  clauses, $i=1,..., 4$, then there
 must be $i$ distinct pairs $\{y_j,z_j\}$ of clones of variables
 in $K_{n,q}$
 such that $\ell(y_j, z_j) < \infty$ and all the $\ell(y_j, z_j)$
 intermediate clauses of each pair are in $F(n,p)$, $j=1,..,i$.
  The probability
 of such event is at most
$$  O\bb( { M_K \choose 8} \sum_{\ell_1, ..., \ell_i\geq 1}
(p/q)^{\ell_1+\cdots + \ell_i} \pr^* [\ell(y_j,z_j) =\ell_j,
~j=1,...,i|SIM] \bb)
$$
for fixed $i$ distinct pairs $\{y_i, z_i\}$ of clones of variables
in $K_{n,p}$, $i=1,..,4$. Since    $M_K =O(\log^3 (\gs^3 n))$ and
  \bean P_{i}&:= &  \sum_{\ell_1, ...,
\ell_i\geq 1} (p/q)^{\ell_1+\cdots + \ell_i} \prc [\ell(y_j,z_j)
=\ell_j, ~j=1,...,i] \\
&\leq &   \sum_{ \ell_1, ..., \ell_i\geq 1}^{
\frac{10}{\gs}\log(\gs^ 3n)} (1-\gs)^{\ell_1+\cdots + \ell_i} \prc
[\ell(y_j,z_j) =\ell_j,
~j=1,...,i] + (1-\gs)^{\frac{10}{\gs}\log(\gs^ 3n)} \\
&  \leq&   \bb(\frac{ 1+o(1)}{8\gs n^{1/3} \log^2 (\gs^3 n)}\bb)^i
+ (\gs^3 n)^{-10}
 , \eean
 that
  $$ \pr [ |F_K ({n,p}) | \geq 4] \leq  (\gs^3 n)^{-4/3+o(1)}. $$

Thus, it is enough to estimate the probability of $|F_K(n,p)|=2,3$
since kernels must have two or more clauses. For two variables
$w,x$ in $K_{n,q}$,
let  $A_{wx}$ be the event of $K_{n,p}=\{ w,x\}$, $|F_K(n,p)|=3$
and  no  variables in $K_{n,q}$ are in  intermediate clauses, and
let $B_{wx} $ be the event that, in addition to $A_{wx}$, each of
the three clauses in $F_{K} (n,p)$ has at least $1$ but not more
than $\frac{10}{\gs}\log (\gs^3 n)$ intermediate clauses. Clearly,
for $K_{n,q} (1,2) := C_{n,q} (1,2)\cup C_{n,q} (2,1)$,
  \beq{nlower} \pr[ |K_{n,p} |=2, |F_K (n,p)
|=3] \geq \pr \Big[ \bigcup_{w,x\in K_{n,q} (1,2)\atop w\not= x}
B_{wx} \Big]~.
 \eeq
For an upper bound, if $|K_{n,p} |=2, |F_K (n,p) |=3$ but
$\cup_{w,x\in K_{n,q} \atop w\not= x} A_{wx}$ does not occur, then
 at least $4$ distinct pairs of clones of variables in
$K_{n,q}$ must have finite length and all corresponding
intermediate clauses of them must be in $F(n,p)$. This probability
is at most $O({M_K \choose 8} P_4) = (\gs^3 n)^{-4/3+o(1)}. $ The
probability of   $\cup_{\{w, x\} \not\sub K_{n,q} (1,2)} A_{wx} $
 may be bounded by
$$\pr \bb[ \bigcup_{\{w, x\} \not\sub K_{n,q} (1,2)} A_{wx}\bb] =
O(|K_{n,q}| (|K_{n,q}|-|K_{n,q} (1,2)|) P_3)= o(\log^6 (\gs^3 n)
P_3)= o((\gs^3 n)^{-1}),$$ for  $|K_{n,q}|-|K_{n,q} (1,2)|
=o(\log^3 (\gs^3 n)). $ Finally, $$\sum_{wx\in K_{n,q} (1,2) \atop
w\not= x } \pr[ A_{wx}\sm B_{wx} ] = O\bb(|K_{n,q}|^2\frac{
P_2}{M_C}\bb) + O\bb(|K_{n,q}|^2(1-\gs)^{\frac{10}{\gs} \log
(\gs^3 n)}\bb) = O( (\gs^2 n)^{-1} ).$$ All together, we have   $$
\pr \bb[ \, |K_{n,p}| =2, |F_K (n,p)|=3] = \pr \Big[
\bigcup_{wx\in K_{n,q}(1,2) \atop w\not=x } B_{wx} \Big]+ o((\gs^3
n)^{-1}).
 $$


Moreover,  as $$ 0\leq  \sum_{wx\in K_{n,q}(1,2)\atop w\not=x }
\pr [B_{wx} ]-\pr \Big[ \bigcup_{w,x\in K_{n,q}(1,2) \atop w\not=x
} B_{wx} \Big] \leq \sum_{w,x,w',x'\in K_{n,q}(1,2) \atop w\not=x,
w'\not=x', \{w,x\}\not= \{w', x'\}} \pr [ B_{wx}\cap B_{w'x'} ]$$
and
$$
\sum_{w,x,w',x'\in K_{n,q}(1,2) \atop w\not=x, w'\not=x',
\{w,x\}\not= \{w', x'\}} \pr [ B_{wx}\cap B_{w'x'} ] = O(
|K_{n,q}|^4 P_4)=o((\gs^3 n)^{-1}), $$ we deduce that
$$ \pr \bb[ \, |K_{n,p}|
=2, |F_K (n,p)|=3 \bb] = \sum_{wx\in K_{n,q}(1,2) \atop w\not=x }
\pr[B_{wx} ]+ o((\gs^3 n)^{-1}).
 $$

To estimate $\pr[B_{wx} ]$, we  may assume that both of $w$ and
$x$ are of type $(2,1)$, after exchanging the roles  $w,x$ with
$\bar{w}, \bar{x}$ if needed. In the random configuration, there
are $5\cdot3\cdot 1$ ways to match the $6$ clones of $w$ and $x$.
In each case $i=1,..., 15$, let $B_i (\ell_1, \ell_2, \ell_3)$ be
the event that there are $\ell_1, \ell_2, \ell_3$ intermediate
variables in $K_{n,q}$ for the three matches. Then
$$ \pr[ B_{wx}] = \sum_{i=1}^{15}
\sum_{\ell_1, \ell_2, \ell_3\geq 2}^{\frac{10 }{\gs}\log (\gs^3
n)}(1-\gs-o(\gs))^{\ell_1+\ell_2+\ell_3}\prc [B_i (\ell_1, \ell_2,
\ell_3)|SIM].$$ For $\ell_1, \ell_2, \ell_3$ in the above range,
\raf{esim} and \raf{nolarge}   give
  \bean
\prc [B_i (\ell_1, \ell_2, \ell_3)|SIM] &=& \frac{\prc[B_i
(\ell_1, \ell_2, \ell_3)\cap SIM]}{\prc[SIM]} \\
 & = &
 \bb(\frac{ 1+o(1)}{8n^{1/3}
\log^2 (\gs^3 n)}\bb)^{\! 3}\, \,  \frac{\prc[SIM']}{\prc[SIM]}\\
&=& \bb(\frac{
1+o(1)}{8n^{1/3} \log^2 (\gs^3 n)}\bb)^{\! 3},
 \eean where $SIM'$ is the event that the random configuration on the
 $M_C - 6 - 2(\ell_{1}+ \ell_2+\ell_3)$ clones is simple. Hence
$$ \pr[ B_{wx}] =  \frac{15+o(1)}{8^3 \gs^3 n \log^6 (\gs^3 n)},$$
and
$$ \pr \bb[ \, |K_{n,p}| =2, |F_K (n,p)|=3 \bb]
= {|K_{n,q} (1,2) |\choose 2} \frac{15+o(1)}{8^3 \gs^3 n \log^6
(\gs^3 n)} + o( (\gs^3 n)^{-1}) = \frac{15+o(1)}{16 \gs^3 n } . $$

If $|K_{n,p}| =1$, then there are at least two clauses in $F_K
(n,p)$. Appealing directly to $F(n,p)$ with $p=
\frac{1-\gs}{2n-1}$,
$$ \pr [ |K_{n,p}| =1 ] \!= \!O\bb(n\! \! \!
\sum_{\ell_1, \ell_2\geq 1 } { n \choose \ell_1-1} 2^{\ell_1-1}
(\ell_1-1)! { n \choose \ell_2-1} 2^{\ell_2-1}(\ell_2-1)!
\bb(\frac{1-\gs}{2n-1}\bb)^{\ell_1+\ell_2}\bb)=O\bb(\frac{1}{\gs^2
n}\bb),$$ where $\ell_1$ and $\ell_2$ represent the numbers of
intermediate clauses.

When the event $B_{wx}$ occurs for variables $w,x$ of type
$(2,1)$, the only way (out of the 15 ways) it directly makes the
formula unsatisfiable is the case that the two clones of each of
$w$ and $x$ are matched and the two clones of $\bar{w}$ and
$\bar{x}$ are matched. Therefore, the same argument yields
$$ \pr [ F(n,p) ~\mbox{is UNSAT}] = \frac{1+o(1)}{16 \gs^3 n }. $$
\qed

 \bn

\old{
\begin{cor} \label{winupper} If $\gl = 1+ \gs$ with $\gs \gg n^{-1/3}$ is
uniformly bounded from above,
$$\pr [\,  F_{PC} (\gl)   ~\mbox{is SAT} \, \, ]
=  e^{-\Theta (\gs^3 n)} . $$
\end{cor}
} 

\renewcommand{\gt}{t}
\noindent {\bf Proof of Theorem \ref{sub2sat}} At $\gL_S$, $F_K
(n,p)$ has  $\Omega (\thl^3 n)$ variables $1-e^{-\Omega(\thl^3 n)}
$ by Theorem \ref{main} and the convention (see Section 1).
Exchanging the roles of $x$ and $\bar{x}$, if necessary, we may
assume that the number of $x$-clones is at least as large as the
number $\bx$-clones for all variables $x\in K_{n,p} $. For the
lower bound, let $Y$ be the set of clones of $x$'s and $Z$ be the
set of clones of $\bx$'s. Then
$$ |Y| \geq |Z| . $$
We now consider  the event that all clones in $Z$ are matched to
clones in $Y$, in which  case, $(1,..., 1)$ is a satisfying
assignment. The probability of the event is
$$ \frac{|Y|}{|Y|+|Z|-1}\, \frac{|Y|-1}{|Y|+|Z|-3} \cdots \frac{|Y|-|Z|+1}{|Y|-|Z|+1}
\geq 2^{-|Z|} \geq e^{-O(\thl^3 n)}=e^{-O(\gs^3 n)} . $$

For the upper bound, we may assume  $\gs \leq 0.01$. Since the
probability decreases as $\gs$ increase, once the probability is
at most $e^{-\Omega (\gs^3 n)}$ for $\gs=0.01$, the probability is
at most $e^{-\Omega ( n)}$ for larger $\gs$'s.
 Corollary \ref{corekernel} implies
that, with probability $1-e^{-\Omega(\thl^3 n)}$,
   \beq{econd} |K_{n,p} | \geq
0.99\thl^3 n~~{\rm and }~~ M_{n,p} (2,2)+ M_{n,p}(1,3)+ M_{n,p}
(3,1)  \leq 0.01 \thl^3 n. \eeq

It is enough to show the desired bound in the random configuration
satisfying \raf{econd} as $\pr^* [SIM] $ $= \Omega (1)$. We first
take the following procedure to make the problem simpler. Remove
all the clauses in the kernel $F_K(n,p)$ containing any variable
not in $K_{n,p} (1,2)$ and its negation. (Recall $K_{n,p} (1,2) =
C_{n.p}(1,2) \cup C_{n,p}(2,1)$.) Then \raf{econd} implies that
there are at most $0.01 \thl^3 n$ such clauses. Furthermore, as at
most one variable in $K_{n,p}(1,2)$ is affected by one such
clause, there are at most $0.02 \thl^3 n$ pure clones can be
created. We now apply PLA: Each time a pure clone  is matched, the
number of pure clones decreases by $2$ if it is matched to another
pure clone.
 If it is matched to a non-pure clone,
two clones become pure with probability $1/3$, and  a variable
becomes of type $(1,1)$ with the other probability. This is so
since all non-pure variables are of type $(1,2)$ or $(2,1)$.
Hence,  after each step, the number of pure clones decreases  by
at least $1/3$ in expectation, and increases by  no more than $1$
at any case. The generalized Chernoff bound  implies that  no pure
clone is left within $0.07 \thl^3 n$ steps, with probability
$1-e^{-\Omega(\thl^3 n)}$. Therefore, there are at least
$0.9\thl^3 n$ variables of type $(1,2)$ or $(2,1)$ remain after
PLA stops, with probability $1-e^{-\Omega(\thl^3 n)}$. The desired
bound may be obtained  from the next lemma.

\begin{lem}
Let $F(b)$ be the (multi)formula yielded by the random
configuration  on  $b$ variables of type $(1,2)$ or $(2,1)$. Then
$$ \pr [ F(b) ~\mbox{is SAT} ] =o(e^{-0.02b}), $$
as $b\ra \infty$.
\end{lem}

\mn \pf After exchanging the roles of $x$ and $\bar{x}$ as needed,
we may  assume that all $b$ variables are of type $(2,1)$.
 Notice that an assignment for the $b$ variables may be regarded
 as  a  $0,1$ vector of
length $b$. That is, the $i^{\rm th}$ coordinate of it tells the
truth value of the $i^{\rm th}$ variable, say $x_i$. Suppose an
assignment has exactly $\gt b$ $0$'s. Then, there are $2\gt b +
(1-\gt)b $ clones whose truth values  are set to be $0$. These
clones are called negative. The other clones are set to be $1$ and
will be called positive. The assignment is a satisfying assignment
if and only if there is no edge connecting two negative clones. We
call such an edge {\em bad}. A clause corresponding a bad edge is
also called {\em bad}.

If $F:=F(b)$ is satisfiable, then there are
 assignments that yield  no bad clause. Among those
assignments, we may take one with maximum number of $1$'s. Such an
assignment is called  maximal. Suppose an satisfying assignment
$s=(s_i)$ is maximal. Then, for a variable $x_{_i}$ with $s_i=0$,
the only clone of $\bar{x}_{_i}$, which is  a positive clone (with
respect to $s$), must be connected to a negative clone. Otherwise,
the assignment $s^*$ obtained from $s$ by changing the value of
$s_i$ to $1$ is another satisfying assignment, which implies that
$s$ is not maximal.

 Summarizing, we have the followings. Provided $s$ has $\gt b$ $0$'s,
the number $N$ of negative clones is $2\gt b + (1-\gt)b =(1+\gt)b
$ and the number $M$ of positive clones is $2(1-\gt)b + \gt b
=(2-\gt)b $. If $s$ is a maximal satisfying assignment, then there
is no bad clause and the positive clone of each $\bar{x_i}$ with
$s_{_i}=0$ must be matched to a negative clone. The number $L$ of
positive clones of $\bar{x}_i$ with $s_i=0$ is $\gt b $. Since all
$N$ negative clones must be matched to positive clones, and the
$L$ positive clones mentioned above must be matched to negative
clones, and the number of perfect matchings on $m$ clones for even
$m$ is
$$(m-1)!! = \frac{m!}{2^{m/2} (m/2)!},$$ we have that
\begin{eqnarray*}P(s) &:= & \Pr[ s ~\mbox{is a maximal satisfying assignment}
]\\
&\leq &   \frac{   {M -L \choose N - L} N ! (M -N- 1)!!}{(M+N
-1)!!}
\end{eqnarray*}
 where
$$ N= (1+\gt)b , ~~M=(2-\gt)b
, ~~L =\gt b ,$$ (provided $s$ has $\gt b$ $0$'s). Using Stirling
formula, we have
  \begin{eqnarray*}  P(s) &\leq &  b \exp \Big(
 2(1-\gt) b H(\mbox{$\frac{1}{2(1-\gt)}$}
 )
+ N \ln N + \frac{M-N}{2 } \ln (M-N) - \frac{M+N}{2} \ln (M+N)
 \Big) \\ &=&  b \exp\Big( 2(1-\gt)bH(\mbox{$\frac{1}{2(1-\gt)}$}) + (1+\gt)
b \ln \frac{1+\gt}{3}+\frac{(1-2\gt)b }{2}\ln \frac{1-2\gt}{3}
\Big).
\end{eqnarray*}

Finally, by  ${b \choose \gt b} \leq e^{b H(\gt)}$ and
$$ \max_{0\leq t\leq 1/2} \bb\{ H(\gt) +
2(1-\gt) H(\mbox{$\frac{1}{2(1-\gt)}$}) + (1+\gt)  \ln
\frac{1+\gt}{3} + \frac{1-2\gt }{2}\ln \frac{1-2\gt}{3}\bb\} <
-0.02 ,$$ we have
  $$\pr [ F ~\mbox{is SAT} ]\leq \sum_{t=0}^b
  {b \choose tb} P(s) =
  o( e^{-0.02 b}). $$
  \qed

\include{dstart}



\dvi

\newcommand{\CJM}{Canad. J. Math}
\newcommand{\Combi}{Combinatorica}
\newcommand{\CPC}{Combinatorics, Probability \& Computing }
\newcommand{\DCG}{Discrete Comput. Geom.}
\newcommand{\DM}{Discrete Mathematics}
\newcommand{\DAM}{Discrete Applied Mathematics}
\newcommand{\DTR}{DIMACS Technical Report}
\newcommand{\EJC}{Europ. J. Combinatorics}
\newcommand{\GAC}{Graphs and Combinatorics}
\newcommand{\JCISS}{J.  Combi. Inform. Sys. Sci.}
\newcommand{\PFICEPFS}{Proc. 1st Int'l Conference on Extremal Problems for Finite Sets, {\rm Visegr\'ad, June}}
\newcommand{\JCSS}{J. of Comp.  Sys. Sci.}
\newcommand{\JCT}{J. of Combinatorial Th.}
\newcommand{\JCTA}{J. of Combinatorial Th. (A)}
\newcommand{\JCTB}{J. of Combinatorial Th. (B)}
\newcommand{\JGT}{J. of Graph Th.}
\newcommand{\JTP}{J. of Theoretical Probability}
\newcommand{\MZ}{Mathematische Zeitschrift}
\newcommand{\RSA}{Random Structures and Algorithms}
\newcommand{\SJDM}{Siam J. of Discrete Mathematics}

{\end{document}}
\begin{thebibliography}{99}

\bibitem{A1}
D. Achlioptas. Setting 2 variables at a time yields a new lower
bound for random 3-{SAT} (extended abstract), \textit{Proc.\ 32nd
ACM Symposium on Theory of Computing}, 28--37 (2000).

\bibitem{Dsur} D. Achlioptas. Lower Bounds for Random 3-SAT via
Differential Equations, \textit{Theoretical Computer Science},
{\bf 265}, 159-185 (2001).

\bibitem{AS}
D. Achlioptas and G.B. Sorkin. Optimal myopic algorithms for
random 3-SAT, \textit{Proc.\ 41st Symposium on the Foundations of
Computer Science}, 590--600 (2000).

\bibitem{AKKK}  D. Achlioptas, L. Kirousis, E. Kranakis, and
D. Krizanc.  Rigorous Results for (2+p)-SAT, {\em Theoretical
Computer Science}, {\bf 265}, 109--129 (2001).










\bibitem{BBK}
A. B\'ek\'essy, P. B\'ek\'essy and J. Koml\'os,  Asymptotic
enumeration of regular matrices, {\em Studia Sci. Math. Hungar.},
 {\bf 7} (1972), 343--353.


\bibitem{BC}
E. Bender and R. Canfield,  The asymptotic number of labeled
graphs with given degree sequences, {\em  J. Combinatorial Theory
Ser. A} {\bf 24} (1978), 296--307.

\bibitem{Bo}
B. Bollob\'as, A probabilistic proof of an asymptotic formula for
the number of labelled regular graphs, {\em  European J. Combin.}
{\bf 1} (1980),  311--316.

\bibitem{BBCKW}
B. Bollob\'as, C. Borgs, J. Chayes, J. H. Kim, D. Wilson, The
Scaling Window of the 2-SAT Transition, {\em \RSA}, {\bf 18}
(2001), 201--256.


\bibitem{BCM90} E.A. Bender, E.R. Canfield and B.D. McKay.
The asymptotic number of labeled connected graphs with a given
number of vertices and edges, \textit{Rand.\ Struc.\ Alg.\ }
\textbf{1}, 127--169 (1990).

\bibitem{BFU} A. Broder, A. Frieze and E. Upfal. On the
satisfiability and maximum satisfiablity of random 3-CNF formulas,
\textit{Proc.\ 4th ACM-SIAM Symposium on Discrete Algorithms},
322--330 (1993).










\bibitem{CF} M.T. Chao and J. Franco. Probabilistic
analysis of a generalization of the unit-clause literal selection
heuristics for the $k$ satisfiable problem, \textit{Information
Science} \textbf{51}, 289--314 (1990).





\bibitem{CR} V. Chv\'atal and B. Reed. Mick gets some
(the odds are on his side), \textit{Proc.\ 33rd Symposium on the
Foundations of Computer Science}, 620--627 (1992).


\bibitem{Coo71} S.A. Cook. The complexity of theorem-proving
procedures, \textit{Proc.\ 3rd ACM Symposium on Theory of
Computing}, 151--158 (1971).




\bibitem{DB} O. Dubois and Y. Boufkhad. A general
upper bound for the satisfiablity threshold of random $k$-SAT
formulas, \textit{J.\ Algorithms} \textbf{24}, 395--420 (1997).

\bibitem{DBM}
O. Dubois, Y. Boufkhad, and J. Mandler. Typical random 3-{SAT}
formulae and the satisfiability threshold. Research announcement
at ICTP, Sept.\ 1999.  Two-page abstract
  appears in \textit{Proc.\ 11th ACM-SIAM Symposium on
  Discrete Algorithms}, 126--127 (2000).


\bibitem{EF} A. El Maftouhi and W. Fernandez de la Vega.
On random $3$-sat. \textit{Combin.\ Probab.\ Comput.\ }
\textbf{4}, 189--195 (1995).





\bibitem{FB} E. Friedgut, with appendix by J. Bourgain.
Sharp thresholds of graph properties, and the $k$-sat problem,
\textit{J.\ Amer.\ Math.\ Soc.\ } \textbf{12}, 1017--1054 (1999).




\bibitem{FV} W. Fernandez de la Vega. On random 2-SAT
(revised version), preprint (1998).


\bibitem{FP} J. Franco and M. Paull. Probabilistic
analysis of the Davis-Putnam procedure for solving the
satisfiability problem, \textit{Discrete Applied Mathematics}
\textbf{5}, 77--87 (1983).

\bibitem{FS} A. Frieze and S. Suen. Analysis of two simple
heuristics for a random instance of $K$-SAT, \textit{J.\
Algorithms} \textbf{20}, 312--335 (1996).







\bibitem{Go} A. Goerdt. A threshold for unsatisfiability,
\textit{J.\ Computer and System Sciences} \textbf{53}, 469--486
(1996).





\bibitem{Jea}
S. Janson, Y.C. Stamatiou, and M. Vamvakari.
 Bounding the unsatisfiability threshold of random 3-{SAT},
\textit{Rand.\ Struc.\ Alg.\ } \textbf{17}, 103--116 (2000).




\bibitem{KM} A. Kamath, R. Motwani, K. Palem
and P. Spirakis. Tail bounds for occupancy and the satisfiability
threshold conjecture, \textit{Rand.\ Struc.\ Alg.\ } \textbf{7},
59--89 (1995).

\bibitem{PCM1} J.-H. Kim, Poisson clonind model for random graph,
 manuscript, in http://research.microsoft.com/theory/jehkim/.


\bibitem{PCD} J.-H. Kim, Poisson clonind model for random digraph,
in preparation.

\bibitem{Kea} L. Kirousis, E. Kranakis, D. Krizanc, and Y.
Stamatiou,  Approximating the unsatisfiability threshold of random
formulas, \textit{\RSA}, \textbf{12},  27--38 (1998).



\bibitem{Ka} A. Kaporis, L. Kirousis, Y. Stamatiou, M. Vamvakari
and M. Zito. {\em The unsatisfiability threshold revisited},
submitted.




\bibitem{M} M. Mitzenmacher.
 Tight Threshholds for the Pure Literal Rule, SRC
Technical Note 1997-011.


\old{
\bibitem{MZ1} R. Monasson and R. Zecchina. The
entropy of the $K$-satisfiability problem, \textit{Phys.\ Rev.\
Lett.\ } \textbf{76}, 3881 (1996).

\bibitem{MZ2} R. Monasson and R. Zecchina. Statistical
mechanics of the random $K$-SAT model, \textit{Phys.\ Rev.\ E}
\textbf{56}, 1357--1370 (1997).

}





\bibitem{Wo}
N. Wormald. Some Problems in the Enumeration of Labelled Graphs
PhD thesis, University of Newcastle, 1978.


\bibitem{W}
N. Wormald.  Differential equations for random processes and
random graphs, {\em Ann. Appl. Prob.} \textbf{5}, 1217--1235
(1995).



\bibitem{Zi} M. Zito. Randomised techniques in combinatorial algorithms, Ph.D
thesis, Dept. of Comp. Sci., Univ. of Warwick, 1999.






\end{thebibliography}
